\documentclass[12pt]{elsarticle}

\usepackage{amsmath} 
\usepackage{amssymb}
\usepackage{amsthm}
\usepackage{amsfonts}
\usepackage{bm}
\graphicspath{{figs/}}

\newtheorem{thm}{Theorem}
\newtheorem{lem}[thm]{Lemma}
\newproof{pf}{Proof}
\newcommand{\const}{\mathop{\rm const}\nolimits}

\journal{arXiv}

\begin{document}

\begin{frontmatter}

\title{Splitting Schemes for Non-Stationary Problems with a Rational Approximation for Fractional Powers of the Operator}

\author[nsi,uni]{Petr N. Vabishchevich\corref{cor}}
\ead{vabishchevich@gmail.com}
\ead[url]{https://sites.google.com/view/vabishchevich/}

\address[nsi]{Nuclear Safety Institute, Russian Academy of Sciences, Moscow, Russia}
\address[uni]{North-Eastern Federal University, Yakutsk, Russia}

\cortext[cor]{Corresponding author}

\begin{abstract}
Problems of the numerical solution of the Cauchy problem for a first-order differential-operator equation are discussed.
A fundamental feature of the problem under study is that the equation includes a fractional power of the self-adjoint positive operator.
In computational practice, rational approximations of the fractional power operator are widely used in various versions.
The purpose of this work is to construct special approximations in time when the transition to a new level in time
provided a set of standard problems for the operator and not for the fractional power operator.
Stable splitting schemes with  weights parameters are proposed for the additive representation of rational approximation for a fractional power operator.
Possibilities of using similar time approximations for other problems are noted.
The numerical solution of a two-dimensional non-stationary problem with a fractional power of the Laplace operator is also presented.
\end{abstract}



\end{frontmatter}

\section{Introduction} 

Applied mathematical models for describing nonlocal processes characterized by the presence of fractional derivatives are actively discussed
(see, for example, \cite{baleanu2012fractional,uchaikin}).
Among them, we can single out as the most interesting anomalous diffusion models that are associated with nonlocality in space.
The corresponding boundary value problems for equations that include fractional powers of elliptic operators \cite{Pozrikidis16} 
are posed and solved. 
Non-stationary problems with a fractional power operator require special consideration.

In the approximate solution of standard multidimensional problems, a finite element or finite volume approximation of an elliptic operator is used \cite{KnabnerAngermann2003,QuarteroniValli1994} when we have a sparse matrix $A$.
For example, in the spectral definition \cite{birman1987spectral,carracedo2001theory} of the fractional power of the matrix ($A^{\alpha}$ with $\alpha \in (0,1)$), we obtain a full matrix. The critical computational complexity is associated with the calculation of the elements of this and inverse matrices.
The main efforts to build computational algorithms for solving problems with a fractional power operator are associated with the use of approximations $A^{\alpha}$ \cite{higham2008functions}, which are more convenient for implementation.

For an approximate solution of spectral space-fractional diffusion problems, various approaches are used \cite{bonito2018numerical,harizanov2020rev}.
Direct use of the results of the theory of approximation of functions \cite{stahl2003best} is implemented in the works 
\cite{harizanov2020num,harizanov2020anal}, in which the best uniform rational approximation is developed.
Many works (see, for example, \cite{frommer2014efficient,bonito2015numerical,Aceto2019}) are based on the integral representation (Balakrishnan formula \cite{balakrishnan1960fractional}) of an operator's fractional power using specific quadrature formulas. We also used other integral representations \cite{vabishchevich2020}.
The solution to a fractional power operator's problem can be represented as the solution to some auxiliary problems of a larger dimension. 
For example, in articles \cite{nochetto2015pde,nochetto2016pde}, a second-order elliptic problem is used \cite{Caffarelli}, 
and in \cite{vabishchevich2014numerical,duan2018numerical,CegVab2019}, the Cauchy problem for a pseudo-parabolic equation.
In some cases, such approaches, when used in practice, can be interpreted \cite{Hofreither2020} as special variants of rational approximations of a fractional power operator.

In the approximate solution of non-stationary problems for equations with a fractional power operator
\cite{yagi2009abstract}, we focus primarily on the use of unconditionally stable implicit schemes \cite{SamarskiiTheory}.
The problem's features at a new level in time are conveyed by the operator of fractional diffusion and reaction.
Various rational approximations are constructed for these problems (\cite{harizanov2020reac,Aceto2017,Vabishchevich2018}), 
which are repelled by rational approximations of a fractional power operator.
In the article \cite{Vabishchevich2016}, unconditionally stable schemes for problems with a fractional power operator are based on the regularization of explicit schemes. 

This paper considers the Cauchy problem for a first-order evolution equation with a fractional power operator.
We are constructing time approximations based on the direct use of rational approximations of a fractional power operator.
For the operator of the problem, we have an additive representation with simple pairwise permutable operator terms.
Unconditionally stable splitting  two-level schemes are constructed when the transition to a new level in time is provided by solving a set of standard problems.

The paper is organized as follows.
In Section 2, we consider the Cauchy problem for a first-order evolution equation with a fractional degree of a self-adjoint operator in a finite-dimensional Hilbert space.
The problems of computational implementation of conventional two-level schemes are discussed.
We discuss the rational approximations of the fractional power operator in Section 3.
In Section 4, we build splitting schemes and provide stability conditions.
Section 5 is devoted to some more general problems. In Section 6, we formulate a non-stationary boundary value problem for space-fractional diffusion problems.
After finite-difference approximation, we have the Cauchy problem for a first-order evolution equation with a fractional power operator.
We present the results of numerical experiments for a test two-dimensional space-fractional diffusion problem.
The results of the work are summarized in Section 7.

\section{Non-stationary problems with a fractional power operator} 

Let $H$ be a finite-dimensional Hilbert space. The Cauchy problem for the first-order evolution equation 
with a fractional power operator is considered:
\begin{equation}\label{1}
 \frac{d u}{d t} + A^\alpha  u = 0,
 \quad 0 < t \leq T, 
\end{equation} 
\begin{equation}\label{2}
 u(0)= u^0 ,
\end{equation}
when $\alpha \in (0,1)$. 

We will assume that the linear operator $A: H \mapsto H$ is constant (independent of $t$),
self-adjoint and positively definite:
\begin{equation}\label{3}
 \frac{d}{d t} A = A \frac{d}{d t},
 \quad A = A^* \geq \delta I ,
 \quad \delta  > 0,  
\end{equation} 
where $I$ is the identity operator in $H$.
We are looking for a solution $u(t)$ of the equation (\ref{1}) for all $t \in (0,T]$
from $H$ for a given right-hand side of $f(t)$ and
initial condition (\ref{2}).

We arrive at the problem (\ref{1})--(\ref{3}), for example, after discretization with respect to spatial variables in the numerical solution of the initial-boundary value problems of anomalous diffusion.
When using finite-difference approximations \cite{SamarskiiTheory} $u$ there is a grid function defined at the nodes of the computational grid.
Such an example is discussed by us below in Section 6.

The scalar product for $u, v \in H$ is $(u, v)$, and the norm is $\| u \| = (u, u)^{1/2}$.
For a self-adjoint and positive operator $B$, we define a Hilbert space $H_B$ with scalar product and norm $(u, v)_B = (B u, v), \ \| u \|_B = (u, v)_B^{1/2}$.

For an approximate solution of the problem (\ref{1})--(\ref{3}), one or another difference approximation in time is used.
We must inherit the basic properties of the operator-differential problem.
In particular, initial data stability is of crucial importance.

Let us present the simplest a priori estimate for solving the problem (\ref{1})--(\ref{3}), which we will be guided by when considering an approximate solution.
Taking into account the nonnegativity of $A^\alpha$, we have
\begin{equation}\label{4}
 \|u(t)\| \leq \|u^0\|,
 \quad 0 < t \leq T .
\end{equation} 
We should have similar estimates of stability with respect to the initial data for an approximate solution.

We introduce a uniform, for simplicity, grid in time with step $\tau$ and use the notation
$y^n = y(t^n), \ t^n = n \tau$, $n = 0,\dots,N, \ N\tau = T$.
We have problems finding an approximate solution to non-stationary problems with a fractional power operator at a new level in time, even when using the simplest explicit scheme. In this case
\begin{equation}\label{5}
 \frac{y^{n+1} - y^{n}}{\tau } + A^\alpha y^{n} = 0,
 \quad n = 0,\dots,N-1 , 
\end{equation} 
\begin{equation}\label{6}
 y^0 = u^0 .
\end{equation} 
Thus, we must have an efficient computational algorithm for calculating the values of $A^\alpha y^{n}$ at each level in time.
The necessary and sufficient condition for the stability of the scheme (\ref{3}), (\ref{5}), (\ref{6}) has \cite{SamarskiiTheory,SamarskiiMatusVabischevich2002} the form
\[
 \tau \leq \frac{2}{\|A\|^\alpha} .
\] 
Under these conditions, for an approximate solution, there is an a priori estimate
\begin{equation}\label{7}
 \|y^{n+1}\| \leq \|u^{0}\|, 
 \quad n = 0,\ldots, N-1 .
\end{equation} 

Implicit schemes belong to the class of unconditionally stable schemes.
For the problem (\ref{1}), (\ref{2}) we will use a two-level scheme with the weight $\sigma \in (0,1] $, when instead of (\ref{5}) we have the equation
\begin{equation}\label{8}
 \frac{y^{n+1} - y^{n}}{\tau } + A^\alpha (\sigma y^{n+1} + (1-\sigma ) y^{n}) = 0,
 \quad n = 0,\dots,N-1 . 
\end{equation} 
The difference scheme (\ref{6}), (\ref {8}) approximates (\ref{1}), (\ref{2}) with sufficient smoothness of the solution $u(t)$ with first order in $\tau$ for $\sigma \neq 0.5$ and with the second --- for $\sigma = 0.5$ (Crank-Nicolson scheme).
The scheme (\ref{3}), (\ref{6}), (\ref{8}) is unconditionally stable at $\sigma \geq 0.5$.

The solution at a new level in time when applying the scheme  (\ref{6}), (\ref{8}) is found from the solution of the equation 
\[
 (I + \sigma \tau A^\alpha ) y^{n+1} = y^{n} - (1-\sigma )\tau A^\alpha y^{n} .
\] 
Thus, we must first calculate the right-hand side, which includes the problematic term with $A^\alpha y^{n}$, and then (this is the main difficulty) solve the equation with the operator $I + \sigma \tau A^\alpha$. 

The computational implementation of the scheme (\ref{6}), (\ref{8}) for $\sigma > 0$ can be simplified.
We introduce a new value
\[
 y^{n+\sigma } = \sigma y^{n+1} + (1-\sigma ) y^{n} ,
\] 
then equation (\ref{8}) can be conveniently written in the form
\[
 \frac{y^{n+\sigma} - y^{n}}{\sigma \tau} + A^\alpha y^{n+\sigma} = 0,
 \quad n = 0,\dots,N-1 .  
\] 
Because of this, we can first find  $y^{n+\sigma}$ from
\[
 (I + \sigma \tau A^\alpha ) y^{n+\sigma} = y^{n},
\] 
and then
\[
  y^{n+1} = \frac{1}{\sigma } (y^{n+\sigma } - (1-\sigma ) y^{n}) .
\] 
Therefore, we don't need to calculate $A^\alpha y^{n}$.

At present, computational algorithms for solving stationary problems with a fractional power operator 
based on one or another approximation $A^{-\alpha}$ are well developed. Our natural desire is to use these results in the approximate solution of non-stationary problems directly, without explicitly constructing one or another approximation for $(I + \sigma \tau A^\alpha)^{-1}$.

\section{Rational approximation for fractional powers of the operator} 

We focus on using approximations
\begin{equation}\label{9}
  A^{-\alpha} \approx R_m(A; \alpha) ,
\end{equation} 
when the problem with $R_m (A; \alpha)$ is more computationally acceptable.
At present, in computational practice, the most widely used approaches are a rational approximation.
In this case, for $R_m (A; \alpha)$ the representation 
\begin{equation}\label{10}
 R_m (A; \alpha) = \sum_{i=1}^{m} a_i(\alpha) (b_i(\alpha) I + A)^{-1} ,
\end{equation} 
is used.

The variety of approximation options (\ref{10}) (see for example \cite{harizanov2020num,bonito2015numerical,Aceto2019,vabishchevich2020}) is associated with the choice of coefficients
$a_i(\alpha), b_i(\alpha), \ i = 1, \ldots, m$.
For example, when using the integral representation
\begin{equation}\label{11}
 A^{-\alpha} = \frac{\sin(\alpha \pi)}{\pi} \int_{0}^{\infty } \theta^{-\alpha} (A + \theta I)^{-1} d \theta .
\end{equation} 
the coefficients are associated with the nodes and weights of the used quadrature formula.
The main requirement for the choice of coefficients in the representation (\ref{10}) is that they are positive:
\[
 a_i(\alpha) > 0, 
 \quad b_i(\alpha) > 0, 
 \quad i = 1, \ldots, m .
\] 
Under these conditions, for the approximating operator, we have
\begin{equation}\label{12}
 R_m(A; \alpha) =  R^*_m(A; \alpha) > 0 ,
 \quad R_m(A; \alpha) A =  A R^*_m(A; \alpha) . 
\end{equation} 

When solving non-stationary problems, we pass from the original equation (\ref{1}) to some other one.
Two possibilities can be distinguished for applying approximation (\ref{9}).

In article \cite{Vabishchevich2018a}, instead of (\ref{1}), we use (first opportunity) the equation
\[
 A^{-\alpha } \frac{d u}{d t} +  u = 0,
 \quad 0 < t \leq T .
\]
An approximate solution $v(t)$ is found from the solution of the Cauchy problem
\begin{equation}\label{13}
 R_m(A; \alpha) \frac{d v}{d t} +  v = 0,
 \quad 0 < t \leq T ,
\end{equation} 
\begin{equation}\label{14}
 v(0)= u^0 .
\end{equation}
The computational complexity of solving the problem (\ref{13}), (\ref{14}) is associated with the complicated construction of the operator at the time derivative in equation (\ref{13}).

In this work, following \cite{Vabishchevich2016}, we rewrite (second opportunity) equation (\ref{1}) in the form
\[
 \frac{d u}{d t} + A^{-\beta} A u = 0,
 \quad 0 < t \leq T , 
\]
where $\beta = 1 - \alpha$, $\beta \in (0,1)$.
For $v(t)$ we will use the equation
\begin{equation}\label{15}
 \frac{d v}{d t} + R_m(A; \beta) A v = 0,
 \quad 0 < t \leq T .
\end{equation} 

Under the conditions (\ref{12}) for solving problems (\ref{13}), (\ref{14}) and (\ref{14}), (\ref{15}) a priori estimate
\begin{equation}\label{16}
 \|v(t)\| \leq \|u^0\|,
 \quad 0 < t \leq T ,
\end{equation} 
is valid.
It is similar to the estimate (\ref{4}) for solution of the original problem (\ref{1})--(\ref{3}).

\section{Splitting scheme} 

The transition to a new level in time for the problem (\ref{14}), (\ref{15}) is based on using the equation
\begin{equation}\label{17}
 \frac{d v}{d t} + D v = 0,
 \quad t^n < t \leq t^{n+1} ,
\end{equation} 
wherein
\[
 D = R_m(A; \beta) A, 
 \quad  \beta = 1-\alpha .
\] 
Taking into account (\ref{10}), the operator $ D $ satisfies the additive representation
\begin{equation}\label{18}
 D = \sum_{i=1}^{m} D_i,
 \quad  D_i = a_i(\beta) (b_i(\beta) I + A)^{-1} A,
 \quad i = 1, \ldots, m . 
\end{equation} 
For individual operator terms in (\ref{18}) we have
\begin{equation}\label{19}
 D_i = D_i^* > 0,
 \quad D_i D_j = D_j D_i , 
 \quad i, j = 1, \ldots, m .  
\end{equation} 

We represent the solution of equation (\ref{17}) in the form
\[
 v^{n+1} = S v^n, 
 \quad S = \exp(- \tau D) . 
\] 
For the operator of transition from one level in time to another, taking into account (\ref{18}), (\ref{19}), we obtain the following multiplicative representation
\[
 S = \prod_{i=1}^{m} S_i,
 \quad  S_i = \exp(- \tau D_i) ,
 \quad i = 1, \ldots, m . 
\] 
Because of this, the solution to equation (\ref{17}) can be represented as a solution to the sequence of problems:
\begin{equation}\label{20}
\begin{split}
 \frac{d v_i}{d t} + D_i v_i & = 0,
 \quad t^n < t \leq t^{n+1} , \\
 v_i(t^n) & = \left \{
 \begin{array}{ll}
  v(t^n) ,  & i = 1,  \\
  v_{i-1} (t^{n+1}),  & i = 2, \ldots, m ,  \\
 \end{array}
 \right .  \\
v(t^{n+1}) & = v_m (t^{n+1}) .
\end{split}
\end{equation} 
We emphasize that the system of equations (\ref{19}) gives an exact solution to equation (\ref{17}) at times $t^n, \ n = 1, \ldots, N$. 

Approximation of (\ref{20}) (using the notation $w^{n+i/m} \approx v_i(t^{n+1}), \ i = 1, \ldots, m$) 
by a two-level scheme leads us to the difference scheme of componentwise splitting \cite{Marchuk1990,VabishchevichAdditive} 
\begin{equation}\label{21}
\begin{split}
 \frac{w^{n+i/m} - w^{n+(i-1)/m}}{\tau} & + D_i (\sigma w^{n+i/m} + (1-\sigma) w^{n+(i-1)/m}) = 0, \\
 \quad  i & = 1, \ldots, m,
 \quad n = 0, \ldots, N-1 ,
\end{split}
\end{equation} 
when setting the initial condition
\begin{equation}\label{22}
 w^0 = u^0 .
\end{equation} 

\begin{thm}\label{1}
The additive operator-difference scheme (\ref{18}), (\ref{19}), (\ref{21}), (\ref{22}) 
is unconditionally stable for $\sigma \geq 0.5$ and for an approximate solution the a priori estimate 
\begin{equation}\label{23}
 \|w^{n+1}\| \leq \|u^0\|,
 \quad n = 0, \ldots, N-1 ,
\end{equation} 
is valid.
\end{thm}

\begin{pf}
Taking into account the relations
\[
 w^{n+(i-\sigma)/m}  = \sigma w^{n+i/m} + (1-\sigma) w^{n+(i-1)/m} ,
\] 
\[
 w^{n+(i-\sigma)/m}  = \left (\sigma - \frac{1}{2} \right ) (w^{n+i/m} - w^{n+(i-1)/m}) + \frac{w^{n+i/m} + w^{n+(i-1)/m}}{2} ,
\] 
we multiply equation (\ref{21}) scalarly by $2 \tau w^{n+(i-\sigma)/m}$. For the first term (\ref{21}) we have
\[
\begin{split}
 \left (\frac{w^{n+i/m} - w^{n+(i-1)/m}}{\tau}, w^{n+(i-\sigma)/m} \right ) & =
 (2 \sigma - 1) \|w^{n+i/m} - w^{n+(i-1)/m}\|^2 \\
 & + \|w^{n+i/m}\|^2 + \|w^{n+(i-1)/m}\|^2 .
\end{split}
\]  
Taking into account the positivity of individual positive operators $D_i, \ i = 1 \ldots, m,$ for $\sigma \geq 0.5$, the estimates 
\[
 \|w^{n+i/m} \| \leq \|w^{n+(i-1)/m} \|, 
 \quad  i = 1, \ldots, m , 
\] 
hold. Thus, we obtain the inequality
\[
 \|w^{n+1} \| \leq \|w^{n} \| , 
\]
from which the proved estimate (\ref{23}) follows.
\end{pf}

The estimate (\ref{23}) is consistent with the estimate (\ref{16}) for problem (\ref{1}), (\ref{2}).
The proposed splitting scheme (\ref{21}) approximates the system of equations (\ref{20}) with an error $\mathcal{O} (\tau^2 + (\sigma -0.5) \tau)$ for sufficiently smooth solutions $v_i(t), \ i = 1 \ldots, m$. 
The study of the estimate of the rate of convergence of the approximate solution to the exact one is carried out in the usual way based on the corresponding estimates of stability on the right-hand side \cite{SamarskiiTheory,SamarskiiMatusVabischevich2002}.

When using the splitting scheme, the solution at a new level in time is determined from equations 
\[
 (b_i(\beta) I + (1+a_i(\beta) \sigma \tau) A) w^{n+i/m} = \chi^{n+(i-1)/m} ,
 \quad  i = 1, \ldots, m, 
\] 
for given right-hand sides
\[
 \chi^{n+(i-1)/m} = (b_i(\beta) I + (1-a_i(\beta) (1-\sigma) \tau) A) w^{n+(i-1)/m},
 \quad  i = 1, \ldots, m . 
\] 
The transition to a new level in time is provided by solving $m$ standard problems with the operator $A + c I, \ c = \const > 0$. 

\section{Generalizations} 

The key feature of the operator-differential equation (\ref{17}), (\ref{18}) is connected with the pairwise commutativity (see (\ref{19})) of the operator terms $D_i, \ i = 1 \ldots, m$.
The consequence of this property is that the sequence of solutions to equations (\ref{20}) gives us an exact solution. The difference approximations for individual equations can be constructed independently of each other.
Let us note some possibilities for constructing splitting schemes for more general problems.

An example of an equation more general than (\ref{1}) is
\begin{equation}\label{24}
 B \frac{d u}{d t} + A^\alpha  u = 0,
 \quad 0 < t \leq T, 
\end{equation} 
with a constant positive and self-adjoint operator $B$ that does not commute with $A$.
When using finite element approximation in space, we associate the matrix mass with the operator $B$.

In an approximate solution of the problem (\ref{2}), (\ref{24}), instead of (\ref{17}), we consider the equation
\begin{equation}\label{25}
 B \frac{d v}{d t} + D v = 0,
 \quad t^n < t \leq t^{n+1} .
\end{equation} 
We can go to the equation
\[
 \frac{d \widetilde{v} }{d t} + \widetilde{D} \widetilde{v}  = 0,
 \quad t^n < t \leq t^{n+1} , 
\] 
for the new unknown $\widetilde{v} (t) = B^{1/2} v$ with the operator
\[
 \widetilde{D} = \sum_{i=1}^{m} \widetilde{D}_i,
 \quad  \widetilde{D}_i =  B^{-1/2} D_i  B^{-1/2},
 \quad i = 1, \ldots, m . 
\] 
The self-adjointness and positivity properties of the operators $\widetilde{D}_i, \ i = 1, \ldots, m,$ are preserved, but pairwise permutability does not take place.

Under these more general conditions, we can use different splitting schemes for multicomponent splitting \cite{VabishchevichAdditive}:
component-wise splitting schemes, regularized splitting schemes, vector schemes.
A direct analogue of the scheme (\ref{21}), (\ref{22}) for equation (\ref{25}) will be the splitting scheme when 
\begin{equation}\label{26}
\begin{split}
 B \frac{w^{n+i/m} - w^{n+(i-1)/m}}{\tau} & + D_i (\sigma w^{n+i/m} + (1-\sigma) w^{n+(i-1)/m}) = 0, \\
 \quad  i & = 1, \ldots, m,
 \quad n = 0, \ldots, N-1 .
\end{split}
\end{equation} 
Similarly to Theorem \ref{1}, we formulate stability conditions.

\begin{thm}\label{2}
For $B = B^* > 0$ and $\sigma \geq 0.5 $, the splitting scheme (\ref{18}), (\ref{19}), (\ref{22}), (\ref{26}) is stable and the solution satisfies the estimate
\begin{equation}\label{27}
 \|w^{n+1}\|_B \leq \|u^0\|_B,
 \quad n = 0, \ldots, N-1 . 
\end{equation} 
\end{thm}

Note that the accuracy of this scheme is of the first-order in $\tau$ for all $\sigma$ and depends on the commutators of the operators
$\widetilde{D}_i$ and $\widetilde{D}_j$, $i, j = 1, \ldots, m$. 
An increase in accuracy to the second-order is achieved \cite{VabishchevichAdditive} by choosing $\sigma = 0.5 $ and organizing calculations according to the rule (Fryazinov-Strang symmetrization) 
\[
 D_1 \rightarrow D_2 \rightarrow \cdots \rightarrow D_m \rightarrow D_m \rightarrow  \cdots \rightarrow D_2 \rightarrow D_1 . 
\] 

The computational implementation of the splitting scheme (\ref{26}) is based on solving the equations
\[
 ((b_i(\beta) I + A) B + a_i(\beta) \sigma \tau A) w^{n+i/m} = \chi^{n+(i-1)/m} ,
 \quad  i = 1, \ldots, m, 
\] 
with the corresponding right-hand sides.
The problem at a new level in time can be difficult due to the presence of the term $A B$.
In particular, owing to the operators' non-commutation $A$ and $B$, the operator $A B$ is non-self-adjoint.

To simplify the computational work, we will build a modification of the scheme (\ref{26}).
We will assume that the operator $B$ is positive definite:
\begin{equation}\label{28}
 B = B^* \geq \gamma I,
 \quad \gamma = \const > 0 . 
\end{equation} 
We will construct additive regularized schemes \cite{VabishchevichAdditive}.

As the generating (primary) scheme, we will use the explicit scheme
\[
\begin{split}
 B \frac{w^{n+i/m} - w^{n+(i-1)/m}}{\tau} & + D_i w^{n+(i-1)/m} = 0, \\
 \quad  i & = 1, \ldots, m,
 \quad n = 0, \ldots, N-1 .
\end{split}
\] 
The unconditionally stable scheme is associated with the perturbation of the operators
\[
 D_i \rightarrow R_i,
 \quad R_i = D_i + \mathcal{O}(\tau),
 \quad i = 1, \ldots, m .  
\] 

In a regularized scheme
\begin{equation}\label{29}
\begin{split}
 B \frac{w^{n+i/m} - w^{n+(i-1)/m}}{\tau} & + R_i w^{n+(i-1)/m} = 0, \\
 \quad  i & = 1, \ldots, m,
 \quad n = 0, \ldots, N-1 ,
\end{split} 
\end{equation} 
for the operators $R_i, \ i = 1, \ldots, m,$ put
\begin{equation}\label{30}
 R_i = a_i(\beta) (b_i(\beta) I + A + \sigma \tau a_i(\beta) A)^{-1} A,
 \quad i = 1, \ldots, m , 
\end{equation} 
with some perturbation parameter $\sigma = \const > 0$.

The stability of this scheme is established on the basis of the following auxiliary statement.

\begin{lem}\label{3}
Let in the scheme
\begin{equation}\label{31}
 B \frac{y^{n+1} - y^{n}}{\tau} + Q y^{n} = 0,
 \quad n = 0, \ldots, N-1 ,
\end{equation} 
the operators $B$ and $Q$ are constant, and
\[
 B = B^* > 0, 
 \quad Q = Q^* \geq 0.
\] 
Then at
\begin{equation}\label{32}
 B \geq \frac{\tau }{2} Q ,
\end{equation} 
the scheme (\ref{6}), (\ref{31}) is stable in $H_B$.
\end{lem}

\begin{pf}
Taking into account the relations
\[
\begin{split}
 y^{n} & = \frac{1}{2} (y^{n+1} + y^{n}) - \frac{1}{2} (y^{n+1} - y^{n}) , \\
 y^{n+1} & = \frac{1}{2} (y^{n+1} + y^{n}) + \frac{1}{2} (y^{n+1} - y^{n}) , 
\end{split}
\] 
multiply equation (\ref{31}) scalarly by $2 \tau y^{n+1}$.
Given (\ref{32}), we have
\[
\begin{split}
 \left ( B \frac{y^{n+1} - y^{n}}{\tau} , 2 \tau y^{n+1} \right ) & = 
 \| y^{n+1}\|^2_B - \| y^{n}\|^2_B + (B (y^{n+1}-y^{n}), y^{n+1}-y^{n}) , \\
 (Q y^{n},y^{n+1}) & = \frac{\tau }{2} (Q (y^{n+1}+y^{n}), y^{n+1}+y^{n}) \\
 & -
 \frac{\tau }{2} (Q (y^{n+1}-y^{n}), y^{n+1}-y^{n}) .
\end{split}
\] 
Adding these equalities, we get
\[
\begin{split}
 \| y^{n+1}\|^2_B - \| y^{n}\|^2_B & + \left ( \left (B - \frac{\tau }{2} Q \right ) (y^{n+1}-y^{n}), y^{n+1}-y^{n} \right ) \\
  & + \frac{\tau }{2} (Q (y^{n+1}+y^{n}), y^{n+1}+y^{n}) = 0 .
\end{split}
\] 
Taking into account the nonnegativity of $Q$, under the inequality (\ref{32}), we will come to an estimate
\[
  \| y^{n+1}\|_B \leq  \| y^{n}\|_B ,
\] 
from which the stability of the scheme (\ref{31}) in $H_B$ follows.
\end{pf}

\begin{thm}\label{4}
The splitting scheme (\ref{22}), (\ref{28})--(\ref{30}) is stable
for $2 \gamma \sigma  \geq 1$, and the solution satisfies the estimate (\ref{27}).
\end{thm}

\begin{pf}
Taking into account lemma \ref{3}, it suffices to check the inequality (\ref{32}), which for (\ref{29}), (\ref{30}) takes the form
\[
 B \geq \frac{\tau }{2} a_i(\beta) (b_i(\beta) I + A + \sigma \tau a_i(\beta) A)^{-1} A,
 \quad i = 1, \ldots, m . 
\] 
For the left side, the inequality (\ref{27}) is used, and for the right side --- 
\[
 \frac{\tau }{2} a_i(\beta) (b_i(\beta) I + A + \sigma \tau a_i(\beta) A)^{-1} A < \frac{1}{2\sigma} I,
 \quad i = 1, \ldots, m . 
\] 
Thus, if the parameter $\sigma$ is specified so that $2 \gamma \sigma \geq 1$, the stability of the splitting scheme (\ref{29}), (\ref{30}) is ensured.
\end{pf}

Splitting schemes for a differential-operator equation
\[
 B \frac{d u}{d t} + A^\alpha  u + C u = f(t),
 \quad 0 < t \leq T, 
\]
with the operator $C \geq 0$ and the given right-hand side $f(t)$ are constructed similarly.

\section{Numerical experiments} 

As the test is considered two-dimensional problem in a square
\[
 \Omega = \{ \bm x  \ | \ \bm x = (x_1,x_2), \ 0 < x_k < 1, \ k = 1,2 \} .
\]
We are looking for a solution $u (\bm x, t)$ of the equation with a fractional power of the Laplace operator:
\[
 \frac{\partial u}{\partial t} + (- \triangle)^\alpha  u = 0 ,
 \quad \bm x \in \Omega ,
 \quad 0 < t \leq T .  
\] 
The boundary and initial conditions are
\[
 u(\bm x, t) = 0,
 \quad \bm x \in \Omega ,
 \quad 0 < t \leq T ,  
\] 
\[
 u(\bm x, 0) = u^0(\bm x),
 \quad \bm x \in \Omega .
\] 

In the domain $\Omega $, we introduce a uniform grid
\[
 \overline{\omega}  = \{ \bm{x} \ | \ \bm{x} =\left(x_1, x_2\right), \quad x_k = i_k h_k, \quad i_k = 0,1,...,N_k,
 \quad N_k h_k = 1, \ k = 1,2 \} ,
\]
where $\overline{\omega} = \omega \cup \partial \omega$ and
$\omega$ is the set of interior nodes, whereas $\partial \omega$ is the set of boundary nodes of the grid.
For grid functions $u(\bm x)$ such that $u(\bm x) = 0, \ \bm x \notin \omega$, we define the Hilbert space
$H=L_2\left(\omega\right)$, where the scalar product and the norm are specified as follows:
\[
\left(u, w\right) =  \sum_{\bm x \in  \omega} u\left(\bm{x}\right)
w\left(\bm{x}\right) h_1 h_2,  \quad 
\| y \| =  \left(y, y\right)^{1/2}.
\]

After approximation in space, we get the problem (\ref{1})--(\ref{3}).  
For $u(\bm x) = 0, \ \bm x \notin \omega$, we take the Laplace grid operator $ A $ in the form
\[
  \begin{split}
  A u = & -
  \frac{1}{h_1^2} (u(x_1+h_1,h_2) - 2 u(\bm{x})- u(x_1-h_1,h_2))  \\ 
  & - \frac{1}{h_2^2} (u(x_1,x_2+h_2) - 2 u(\bm{x}) - u(x_1,x_2-h_2)) , 
  \quad \bm{x} \in \omega . 
 \end{split} 
\] 
For problems with sufficiently smooth solution, 
it approximates the differential operator with the truncation error 
$\mathcal{O} \left(|h|^2\right)$, $|h|^2 = h_1^2+h_2^2$. 
For constant $\delta$ in (\ref{3}) we have (see, for example, \cite{SamarskiiTheory,SamarskiiNikolaev1978})
\[
 \delta  = \sum_{k  =1}^{2}  \frac{4}{h_k ^2} \sin^2 \frac{\pi}{2 N_k } .
\] 

\begin{figure}
\centering
\begin{minipage}{0.32\linewidth}
\centering
\includegraphics[width=\linewidth]{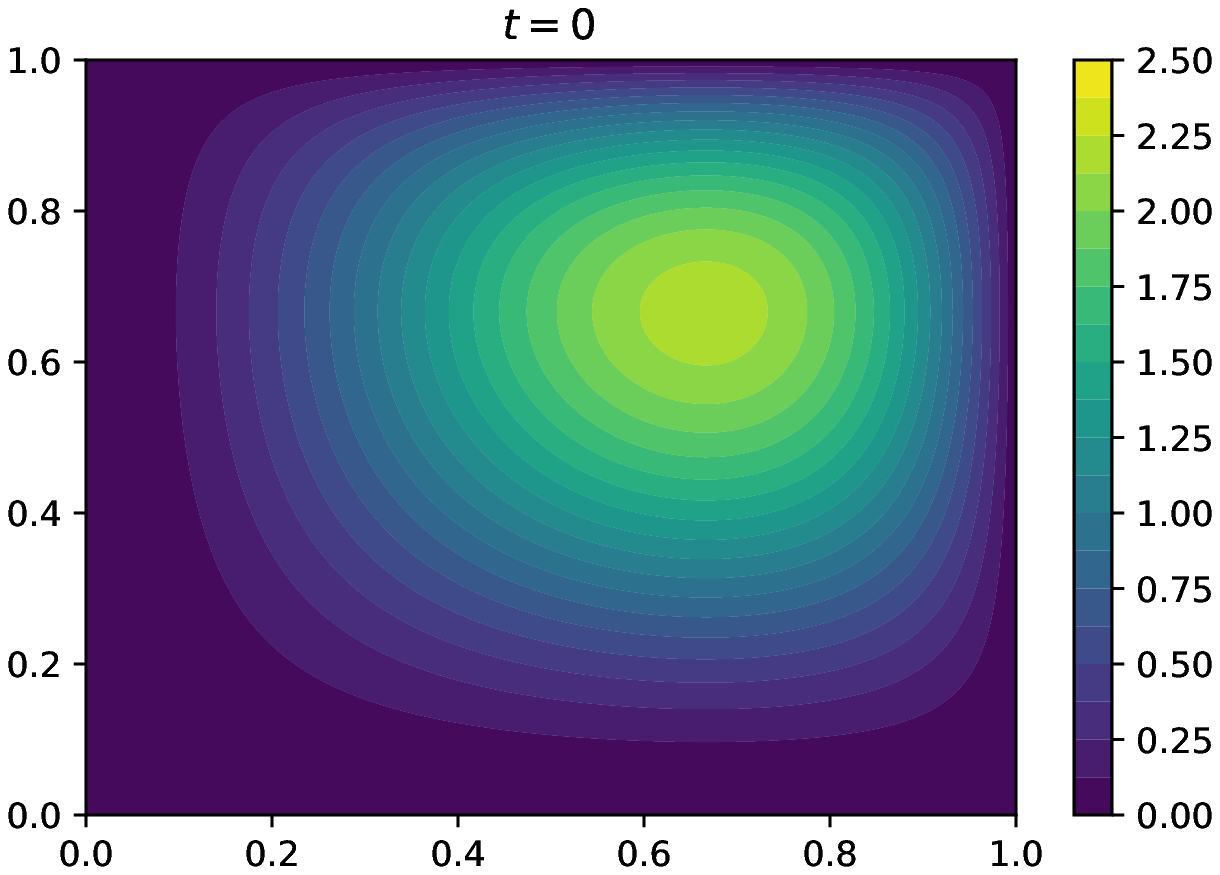}
\end{minipage}
\begin{minipage}{0.32\linewidth}
\centering
\includegraphics[width=\linewidth]{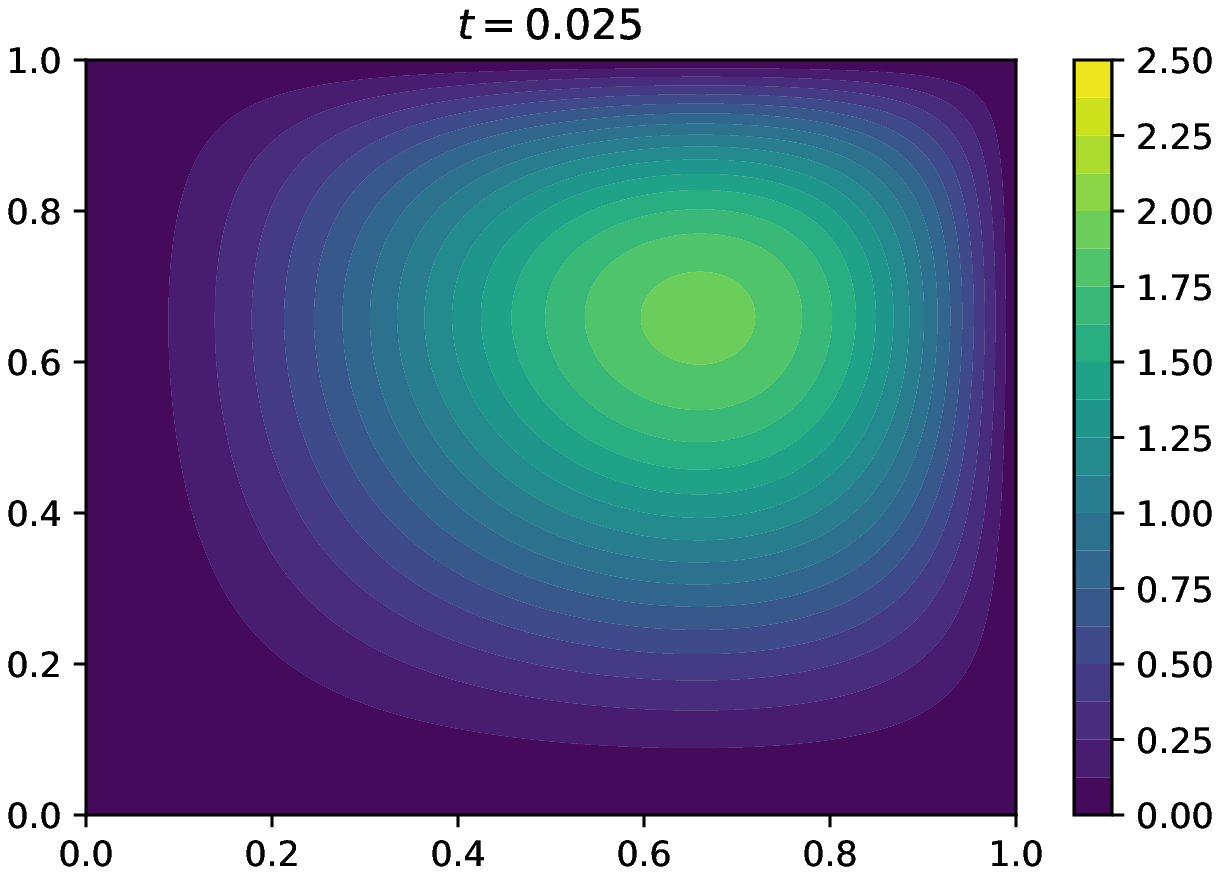}
\end{minipage}
\begin{minipage}{0.32\linewidth}
\centering
\includegraphics[width=\linewidth]{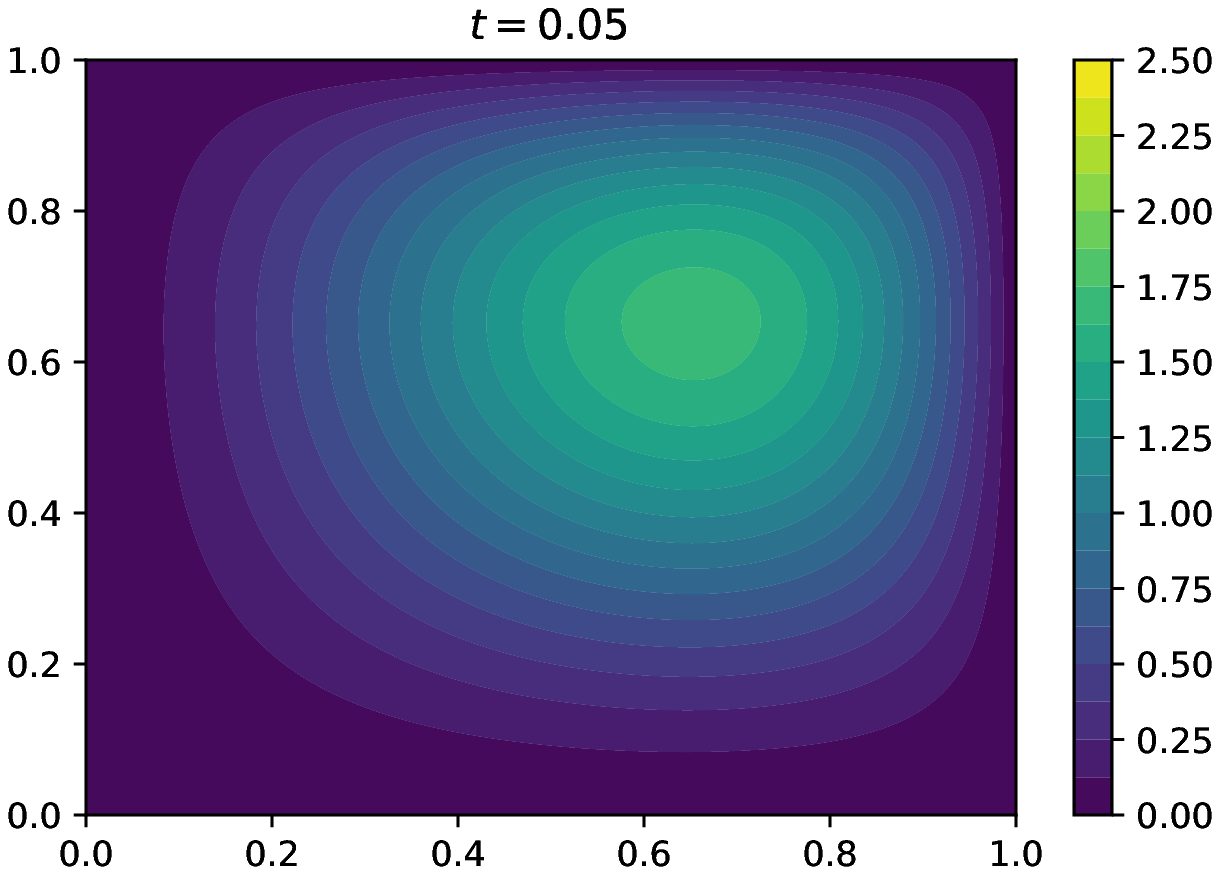}
\end{minipage}
\caption{Exact solution of the problem at separate points in time: $\alpha = 0.5$.}
\label{f-1}
\end{figure}

\begin{figure}
\centering
\begin{minipage}{0.32\linewidth}
\centering
\includegraphics[width=\linewidth]{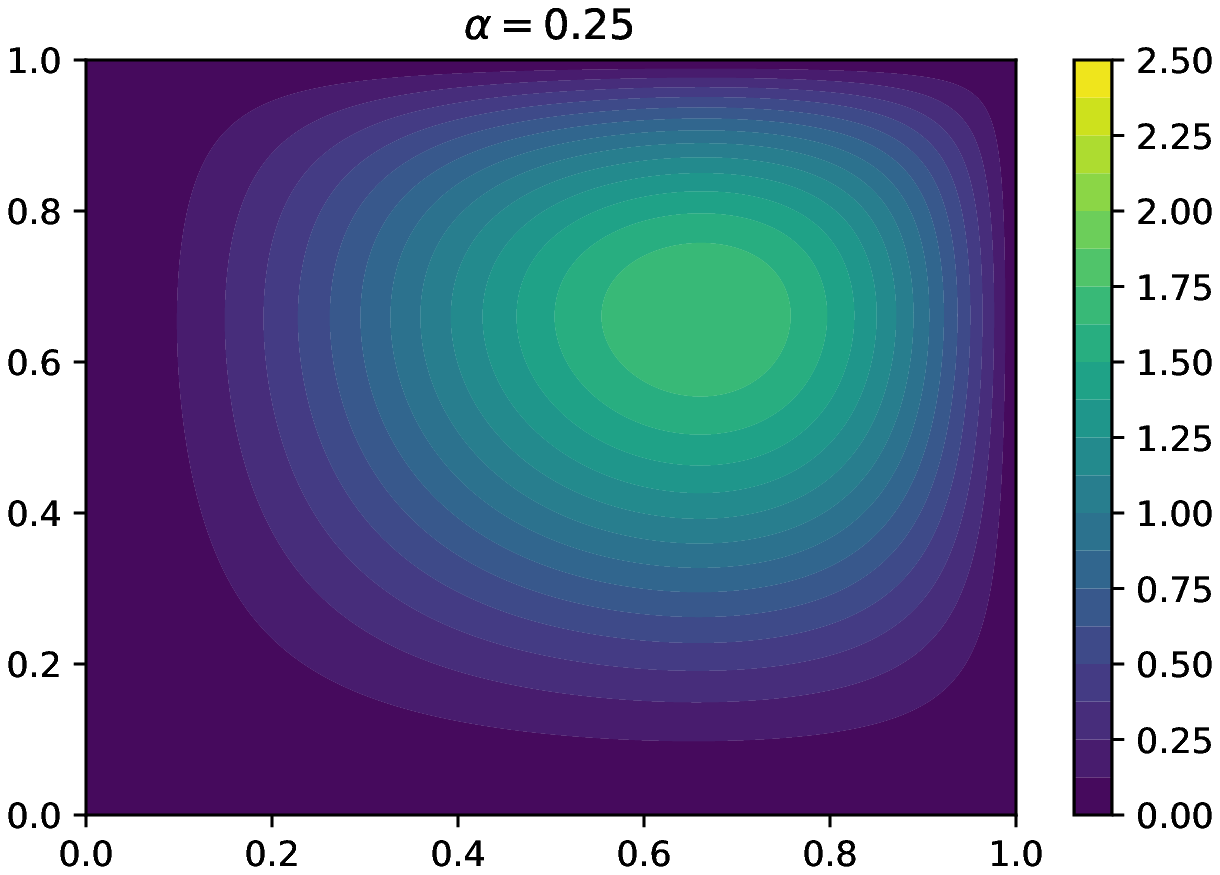}
\end{minipage}
\begin{minipage}{0.32\linewidth}
\centering
\includegraphics[width=\linewidth]{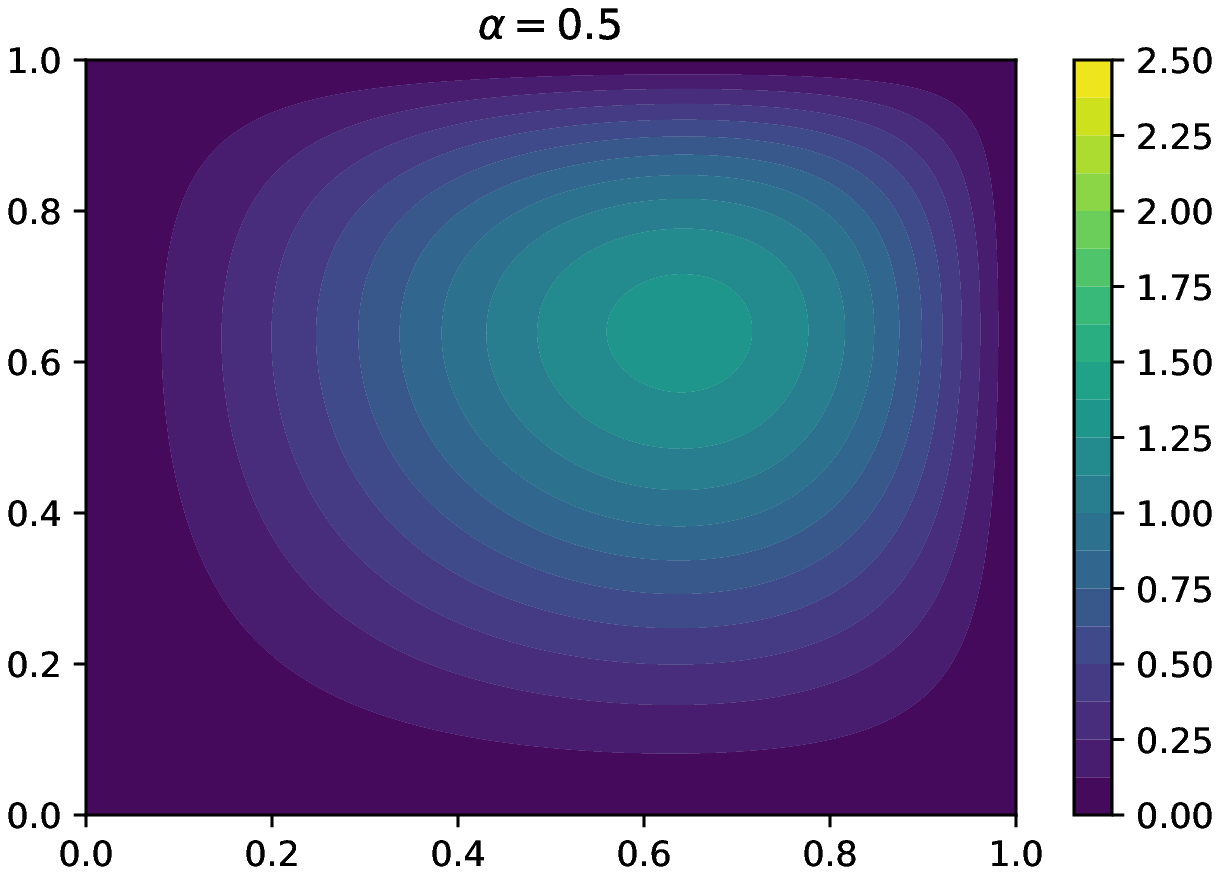}
\end{minipage}
\begin{minipage}{0.32\linewidth}
\centering
\includegraphics[width=\linewidth]{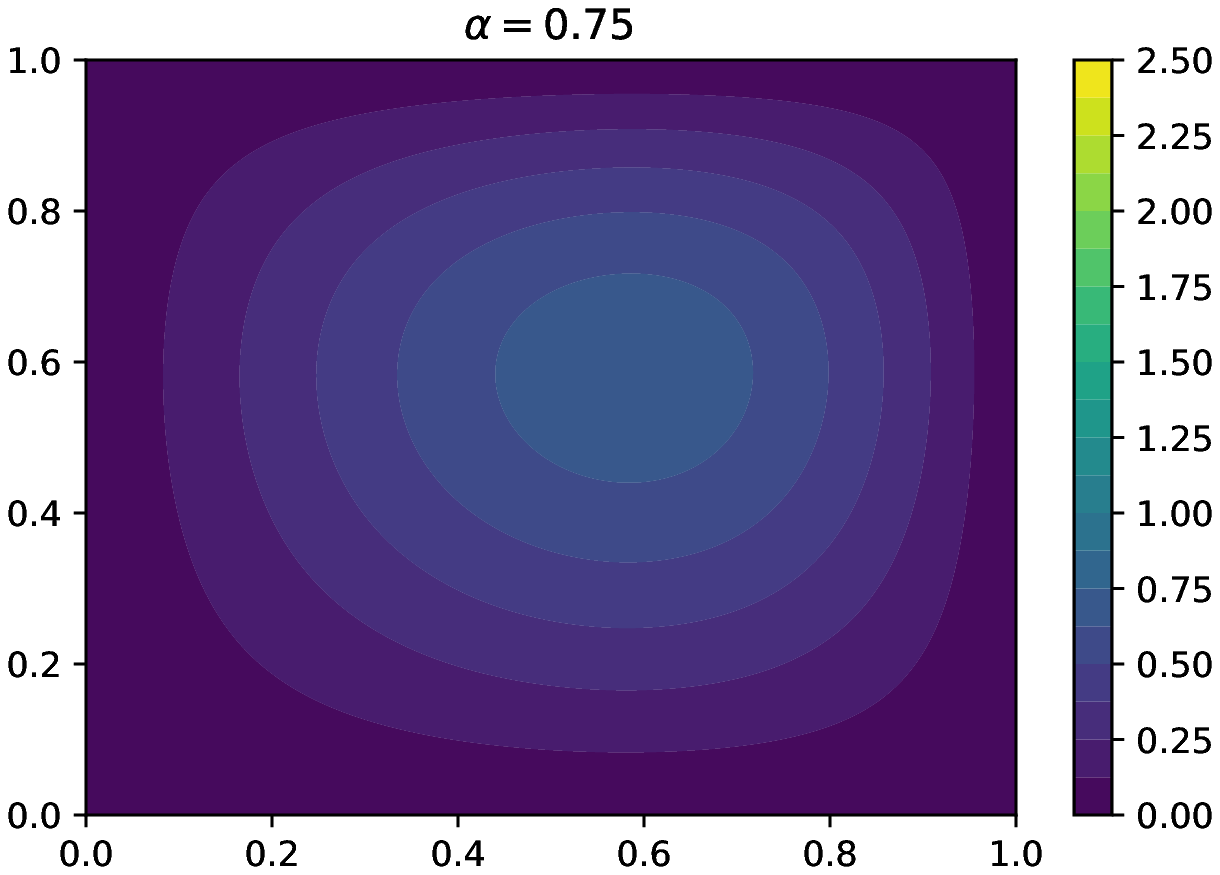}
\end{minipage}
\caption{Exact solution for different values of $\alpha$: $T=0.1$.}
\label{f-2}
\end{figure}

We performed computational experiments for a model problem with the initial condition
\[
 u^0(\bm x) = 100 x_1^2(1-x_1) x_2^2(1-x_2) 
\] 
for $T = 0.1$ and space grid $N_1 = N_2 = 256$. 
The exact solution to the problem (\ref{1}), (\ref{2}) with a value of $\alpha = 0.5$ is shown in Fig.\ref{f-1}.
The influence of $\alpha$ can be seen in Fig.\ref{f-2}.

We can do a rational approximation of $A^{-\beta}$ on the basis of the integral representation (\ref{11}).
We will use Gauss quadrature formulas similarly to \cite{Vabishchevich2018}.
We introduce (see \cite{frommer2014efficient}) a new variable of integration $\eta$ by the relation
\[
 \theta = \mu \frac{1-\eta}{1+\eta},
 \quad \mu > 0.
\] 
From (\ref{11}), we have
\begin{equation}\label{33}
  A^{-\beta } = \frac{2 \mu^{1-\beta} \sin(\pi \beta )}{\pi} \int_{-1}^{1} (1-\eta)^{-\beta}(1+\eta)^{\beta-1}
  \big (\mu (1-\eta) I + (1+\eta) A \big )^{-1} d \eta .
\end{equation}
We use the Gauss-Jacobi quadrature formula  \cite{Rabinowitz} with the weight
$ (1-\eta)^{\tilde{\alpha}}(1+\eta)^{\tilde{\beta}})$:
\[
  \int_{-1}^{1} f(t) (1-\eta)^{\tilde{\alpha} }(1+\eta)^{\tilde{\beta}} d \eta \approx 
  \sum_{i=1}^{m} \upsilon_i f(\eta_i) ,
  \quad \alpha, \beta > - 1. 
\] 
Here $\eta_1, \ldots, \eta_m$ are the roots of the Jacobi polynomial $J_m(\eta; \tilde{\alpha},\tilde{\beta})$
of degree $m$ and $\upsilon_1, \ldots, \upsilon_M$ are weights:
\[
\begin{split}
  \upsilon_i =  & - \frac{2 m + \tilde{\alpha} +\tilde{\beta} + 2}{m + \tilde{\alpha} + \tilde{\beta} + 1} 
  \frac{\Gamma (m + \tilde{\alpha} + 1) \Gamma (m + \tilde{\beta} + 1)}{\Gamma (m + \tilde{\alpha} + \tilde{\beta} + 1) (m+1)!} \\
  & \times  \frac{2^{\tilde{\alpha} + \tilde{\beta}}}{J_m^{\,'}(\eta_i; \tilde{\alpha},\tilde{\beta}) 
  J_{m+1}(\eta_i; \tilde{\alpha},\tilde{\beta})} > 0 ,
\end{split}
\] 
where $\Gamma$ denotes the gamma function.
Thus, for $R_m(A; \beta$), we get
\[
 \tilde{\alpha} = - \beta ,
 \quad \tilde{\beta} = \beta - 1, 
 \quad a_i = \frac{2 \mu^{1-\beta} \sin(\pi \beta )}{\pi} \frac{\upsilon_i}{1+\eta_i} ,
 \quad b_i =  \mu  \frac{1-\eta_i}{1+\eta_i} .
\]
For parameter transformation $\mu \geq \delta$, we set $\mu = 4$.

\begin{figure}
\centering
\begin{minipage}{0.32\linewidth}
\centering
\includegraphics[width=\linewidth]{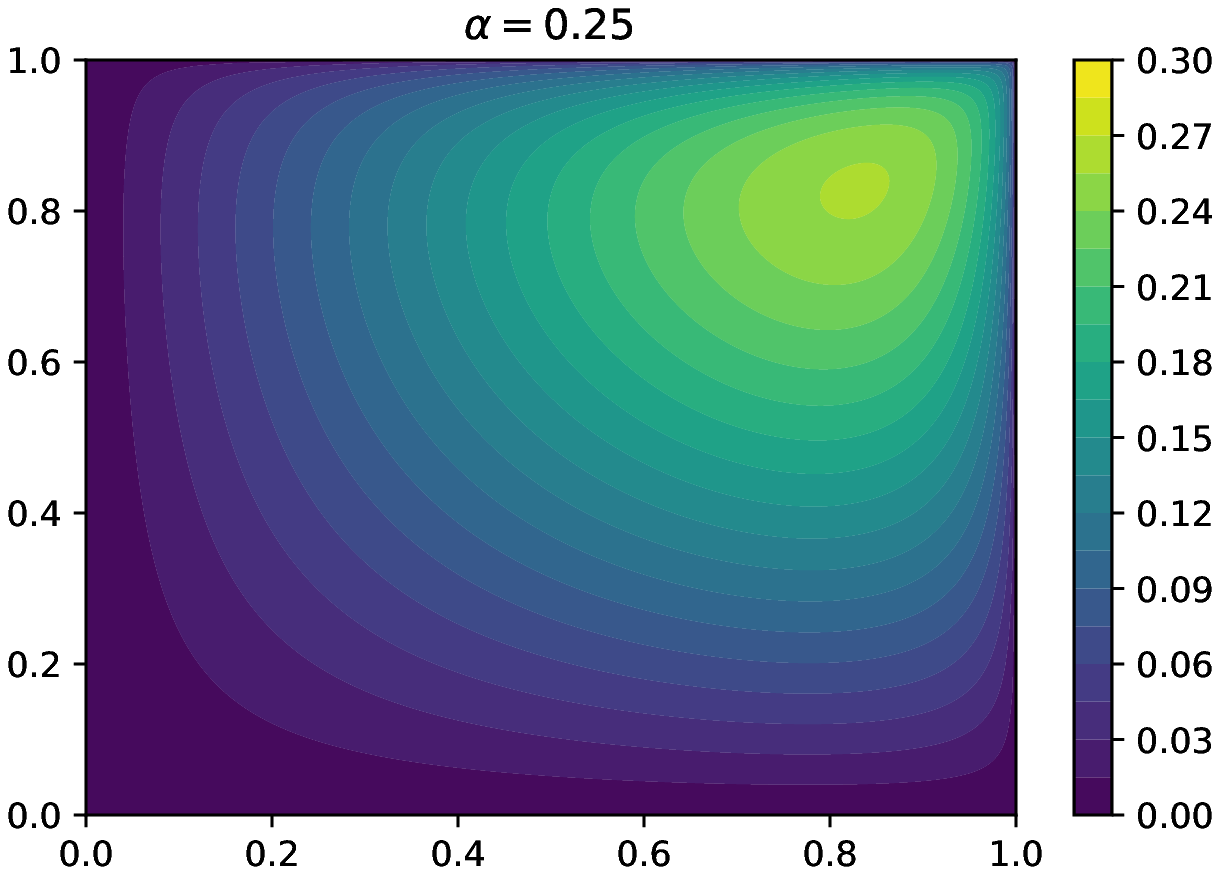}
\end{minipage}
\begin{minipage}{0.32\linewidth}
\centering
\includegraphics[width=\linewidth]{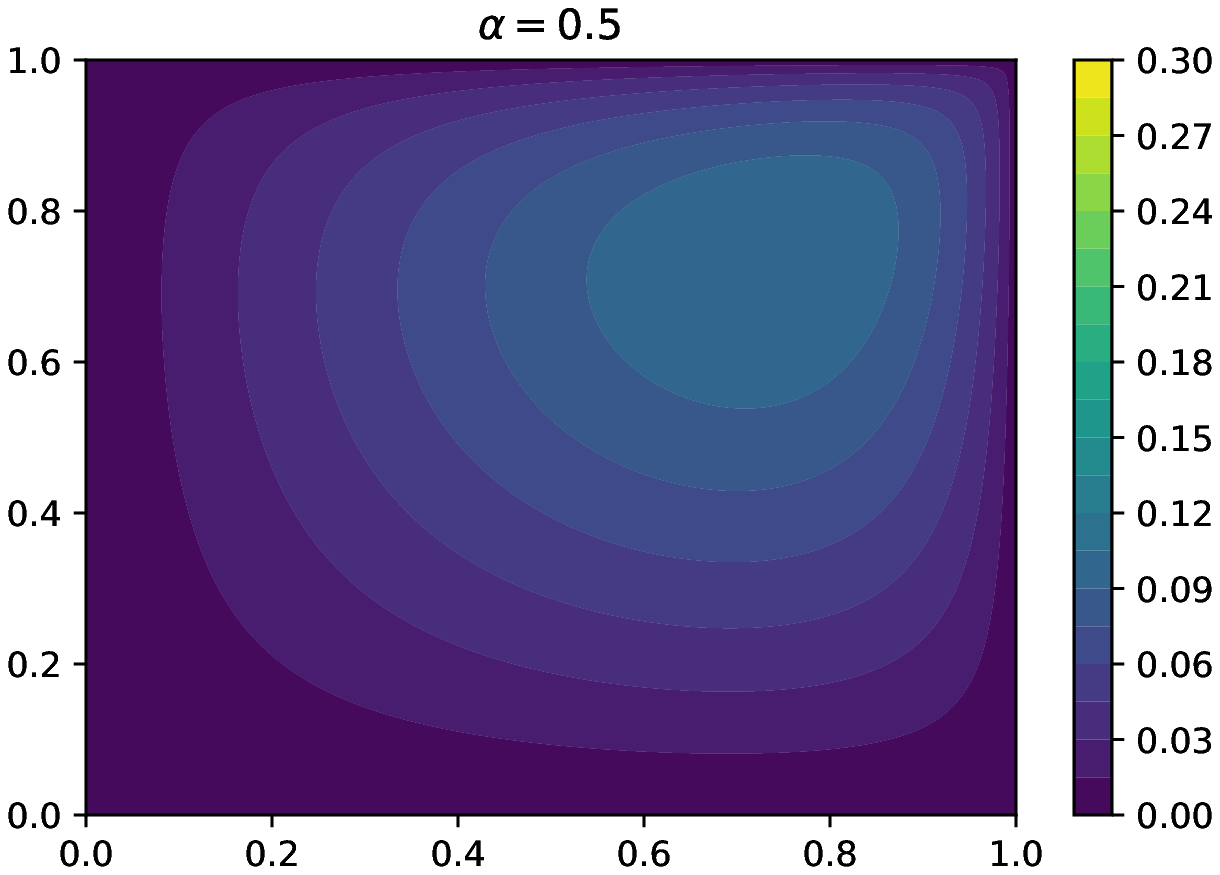}
\end{minipage}
\begin{minipage}{0.32\linewidth}
\centering
\includegraphics[width=\linewidth]{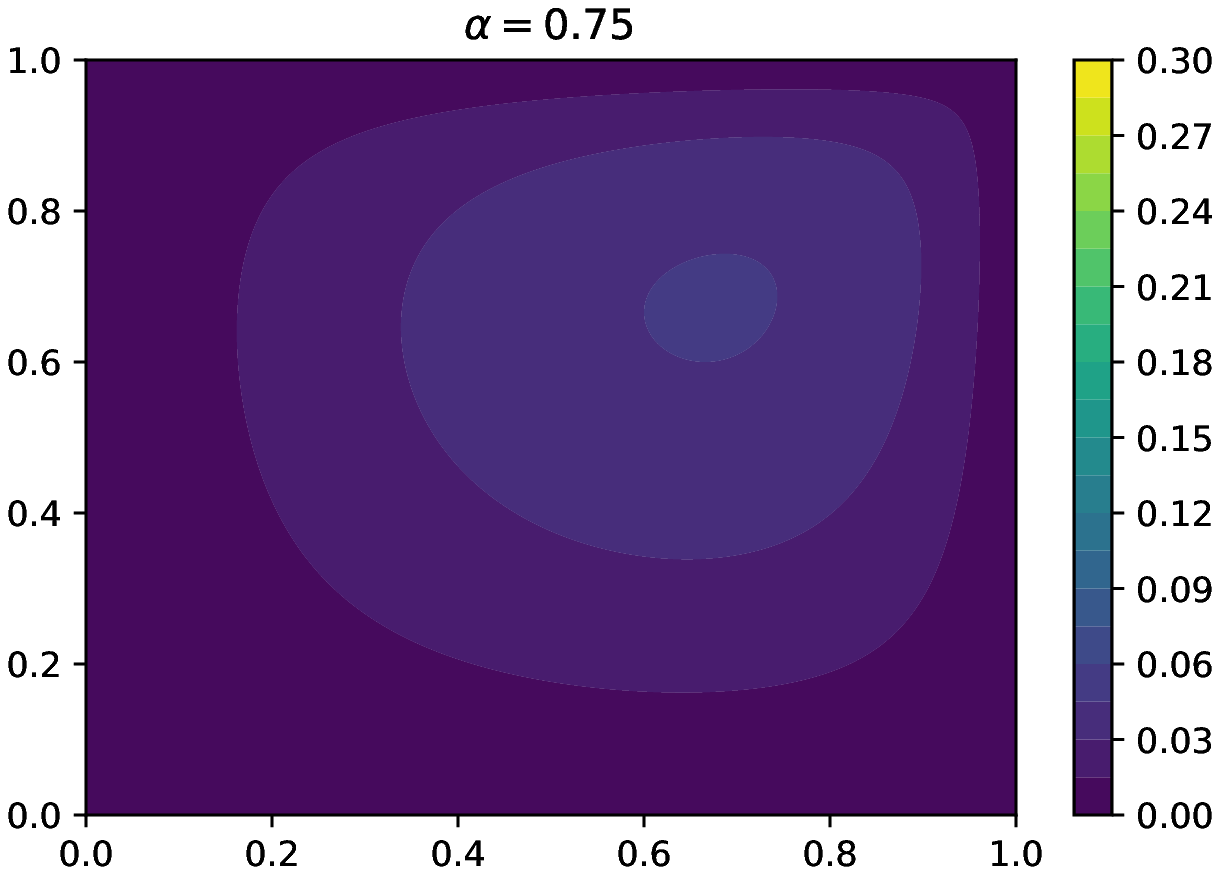}
\end{minipage}
\caption{Exact solution of stationary problem for different values of $\alpha$.}
\label{f-3}
\end{figure}

The accuracy of the rational approximation used will be illustrated by the results of solving the stationary problem
\begin{equation}\label{34}
 A^\alpha v = \varphi , 
\end{equation} 
when the right-hand side is specified on the computational grid $\omega$ by the relation
\[
 \varphi (\bm x) = x_1, x_2,
 \quad \bm x \in \omega . 
\] 
The exact solution for different values of the parameter $\alpha$ is shown in Fig.\ref{f-3}.
Error estimation for an approximate solution $\widetilde{v}$ is performed in $L_2(\omega)$ and $L_\infty (\omega)$:
\[
 \varepsilon_2 = \frac{\|\widetilde{v} (\bm x) - v(\bm x)\|}{ \|v(\bm x)\|},
 \quad  \varepsilon_\infty  = \frac{\|\widetilde{v} (\bm x) - v(\bm x)\|_\infty }{ \|v(\bm x)\|_\infty },
 \quad \|v(\bm x)\|_\infty = \max_{\bm x \in \omega} |v(\bm x) |  .
\] 
We should pay the greatest attention to the dependence of the accuracy of the approximate solution on the number of nodes of the quadrature formula.
As the data from Table~1  show, we can get high accuracy only with a very large number of integration nodes.

\begin{table}[htp]
\label{tab-1}
\caption{Accuracy for the stationary problem: formula (\ref{33}).}
\centering
\begin{tabular}{lllll}
\\
\hline
$\alpha$ &     error                  &  $m = 50$        &    $m = 100$       &    $m = 200$      \\
\hline
0.25     &     $\varepsilon_2     $   &   7.078046e-03   &   1.903738e-03   &   2.231284e-04   \\
         &     $\varepsilon_\infty$   &   9.264885e-02   &   3.370878e-02   &   4.947108e-03   \\
0.5      &     $\varepsilon_2     $   &   1.057186e-03   &   2.550164e-04   &   3.030576e-05   \\
         &     $\varepsilon_\infty$   &   1.142293e-02   &   4.200140e-03   &   6.576673e-04   \\
0.75     &     $\varepsilon_2     $   &   1.047861e-04   &   2.028716e-05   &   2.169714e-06   \\
         &     $\varepsilon_\infty$   &   9.453449e-04   &   3.028348e-04   &   4.497895e-05   \\
\hline
\end{tabular}
\end{table}

Let us illustrate the possibility of a significant increase in accuracy by using other rational approximations of the fractional power of the operator.
The article \cite{vabishchevich2020} uses a new integral representation, the advantages of which are (i) we integrate on a finite interval, (ii) we avoid the singularity of the integrand, (iii) and we control the smoothness of the integrand by choosing the parameter of the integral representation.
Instead of (\ref{33}), we will use the integral representation
\begin{equation}\label{35}
  A^{-\beta } = \frac{\sin(\pi \beta )}{(1-\beta)\pi} \int_{0}^{1} (1-\eta)^{\varkappa -1} 
    \Big (1+ \Big(\varkappa\frac{1-\beta}{\beta} - 1 \Big ) \eta \Big )
    \Big (\eta^{\frac{1}{1-\beta}} I + (1-\eta)^{\frac{\varkappa}{\beta}} A \Big )^{-1} d \eta , 
\end{equation} 
with the parameter $\varkappa > 1$. 

We use Simpson's quadrature formula \cite{Rabinowitz}  to approximate $A^{-\beta }$.
The unit interval $[0,1]$ is divided on $m$ equal parts and the nodes $\eta_i = (i-1)/m, \ i = 1,\ldots, m$ of the quadrature formula are determined, for $\eta = 1$ integrand vanishes.
With this in mind, we define the operators
\[
\begin{split}
 A_i & = d_i(b_i I + A)^{-1} , \\
 d_i & = \frac{\sin(\pi \beta )}{(1-\beta)\pi} (1-\eta_i)^{\varkappa -1 - \frac{\varkappa}{\beta}} 
    \Big (1+ \Big(\varkappa\frac{1-\beta}{\beta} - 1 \Big ) \eta_i \Big ) , \\
 b_i & = \eta_i^{\frac{1}{1-\beta}} (1-\eta_i)^{-\frac{\varkappa}{\beta}} ,
 \quad i = 1,\ldots, m .
\end{split}
\] 
Using Simpson's formula for the coefficients $a_i, \ i = 1,\ldots, m$ in $R_m(A; \beta)$ we get 
\[
 a_1 = \frac{1}{3m} d_1,
 \ a_2 = \frac{4}{3m} d_2, 
 \ a_3 = \frac{2}{3m} d_3, 
 \ \dots,
 \ a_{m-1} = \frac{2}{3m} d_{m-1} ,  
 \ a_m = \frac{4}{3m} d_m .   
\] 
The results of using such a rational approximation with $\varkappa = 5$  for an approximate solution of the problem (\ref{34}) are presented in Table~2.
Comparison with Table~1 shows that the use of rational approximations based on the integral representation (\ref{35}) allows obtaining a solution with much higher accuracy than on the basis of the integral representation (\ref{33}).
We use this variant of the approximating operator $R_m(A; \beta)$ to solve non-stationary problems.

\begin{table}[htp]
\label{tab-2}
\caption{Accuracy for the stationary problem: formula (\ref{35}).}
\centering
\begin{tabular}{llllll}
\\
\hline
$\alpha$ &     error                  &  $m = 50$        &    $m = 100$       &    $m = 200$      \\
\hline
0.25     &     $\varepsilon_2     $   &   3.444934e-07   &   2.380002e-08   &   1.539020e-09   \\
         &     $\varepsilon_\infty$   &   5.568729e-07   &   2.234692e-08   &   1.443448e-09   \\
0.5      &     $\varepsilon_2     $   &   9.686650e-08   &   6.073046e-09   &   3.799713e-10   \\
         &     $\varepsilon_\infty$   &   9.506237e-08   &   5.990690e-09   &   3.748065e-10   \\
0.75     &     $\varepsilon_2     $   &   2.433874e-08   &   1.493598e-09   &   9.359871e-11   \\
         &     $\varepsilon_\infty$   &   5.835720e-08   &   1.491888e-09   &   9.344940e-11   \\
\hline
\end{tabular}
\end{table}

\begin{figure}
\centering
\begin{minipage}{0.45\linewidth}
\centering
\includegraphics[width=\linewidth]{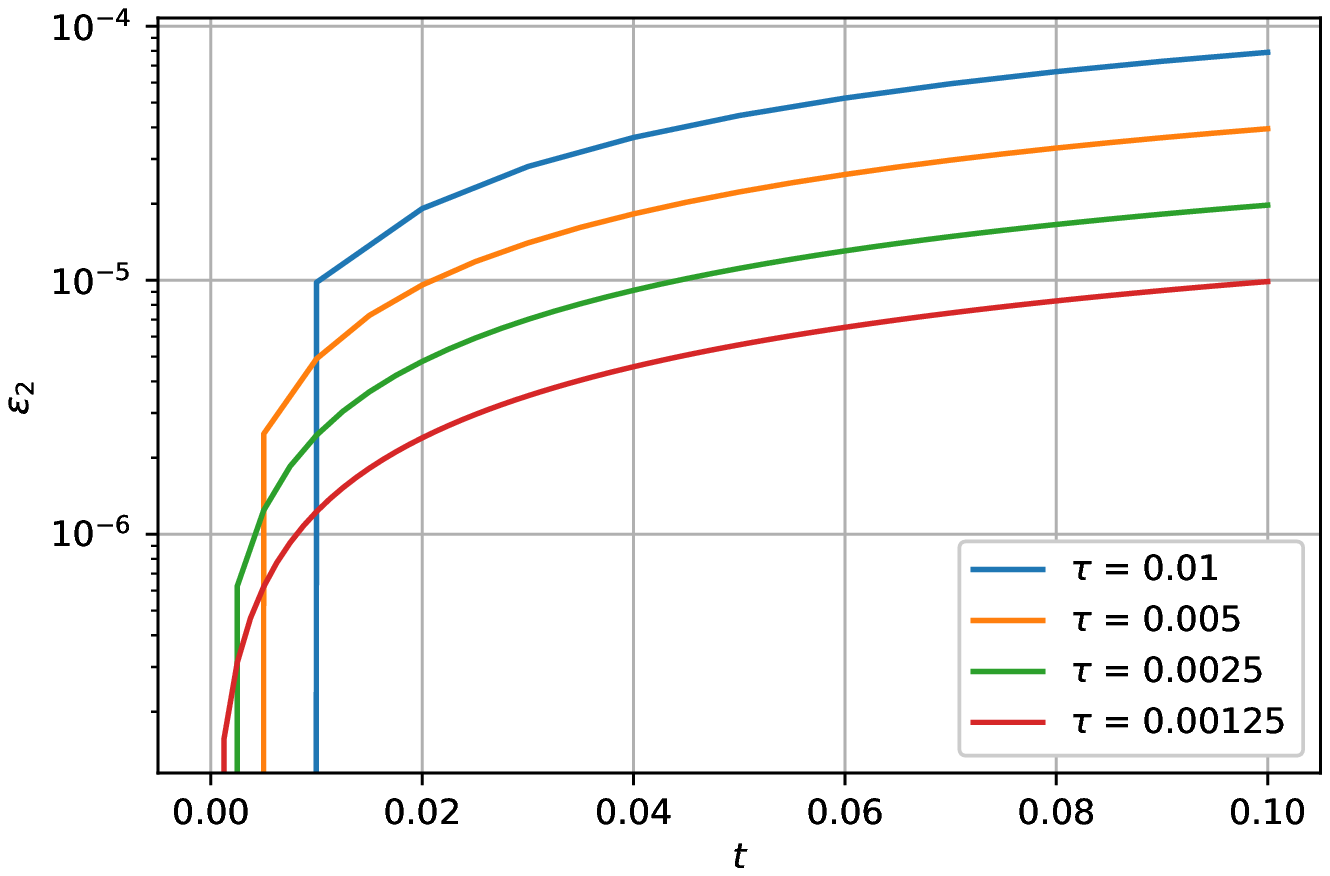}
\end{minipage}
\begin{minipage}{0.45\linewidth}
\centering
\includegraphics[width=\linewidth]{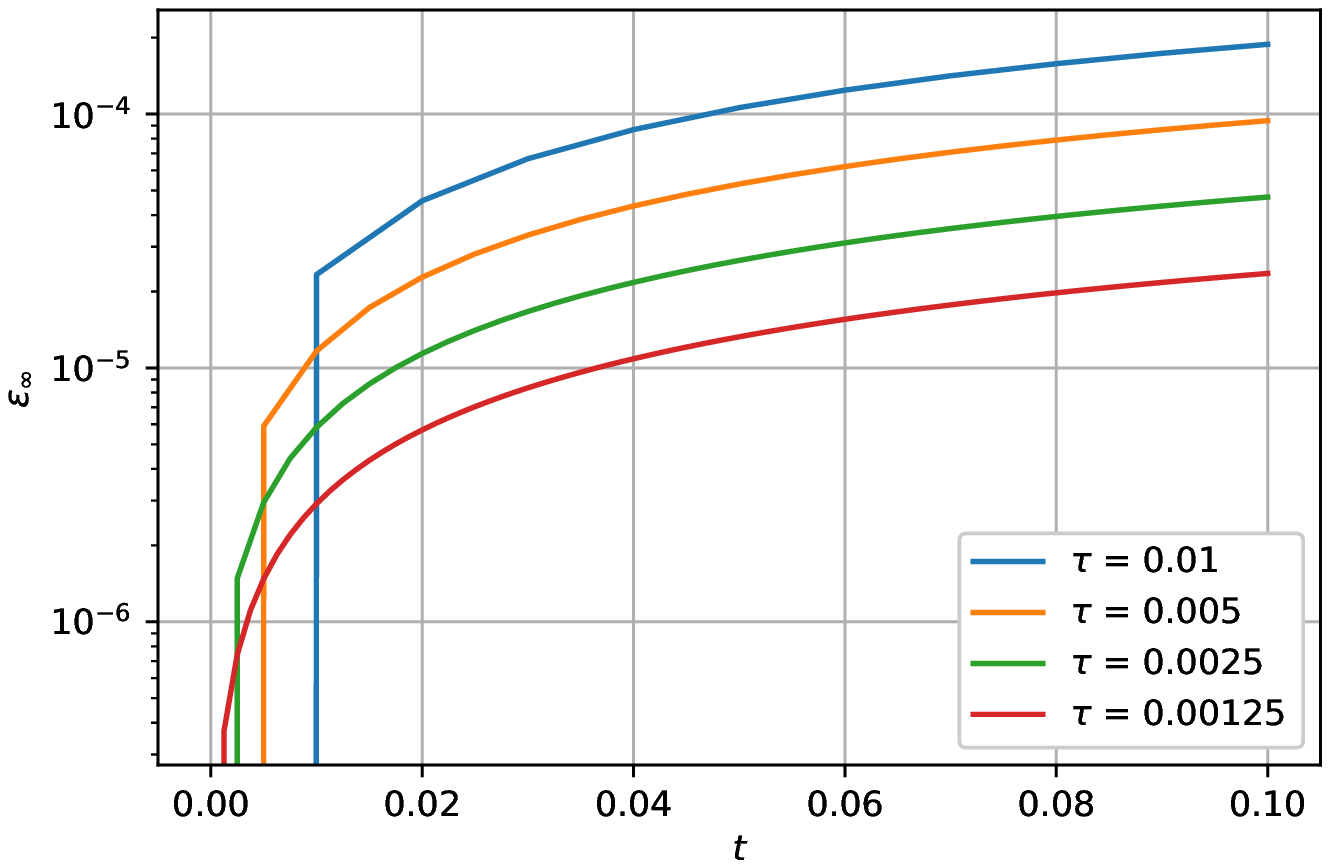}
\end{minipage}
\begin{minipage}{0.45\linewidth}
\centering
\includegraphics[width=\linewidth]{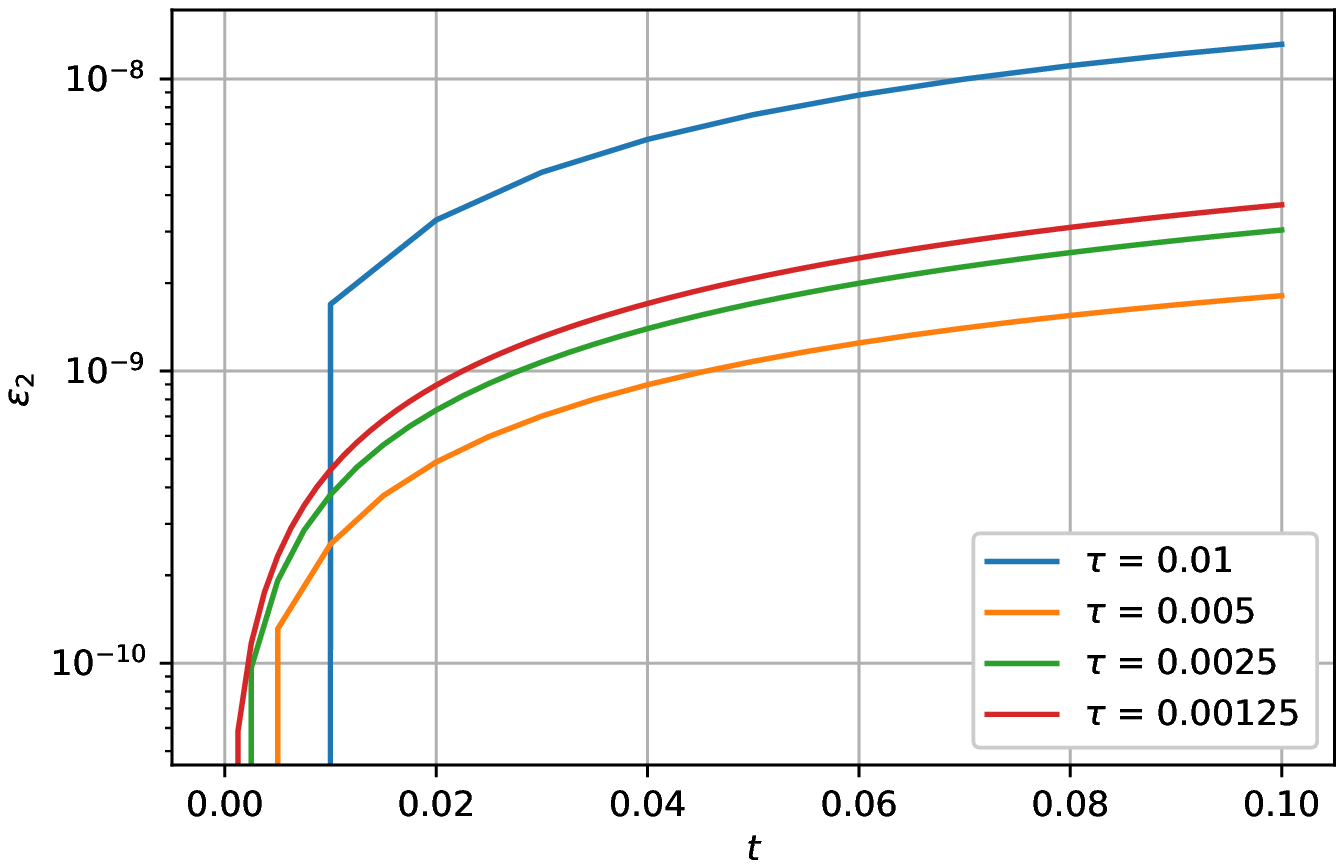}
\end{minipage}
\begin{minipage}{0.45\linewidth}
\centering
\includegraphics[width=\linewidth]{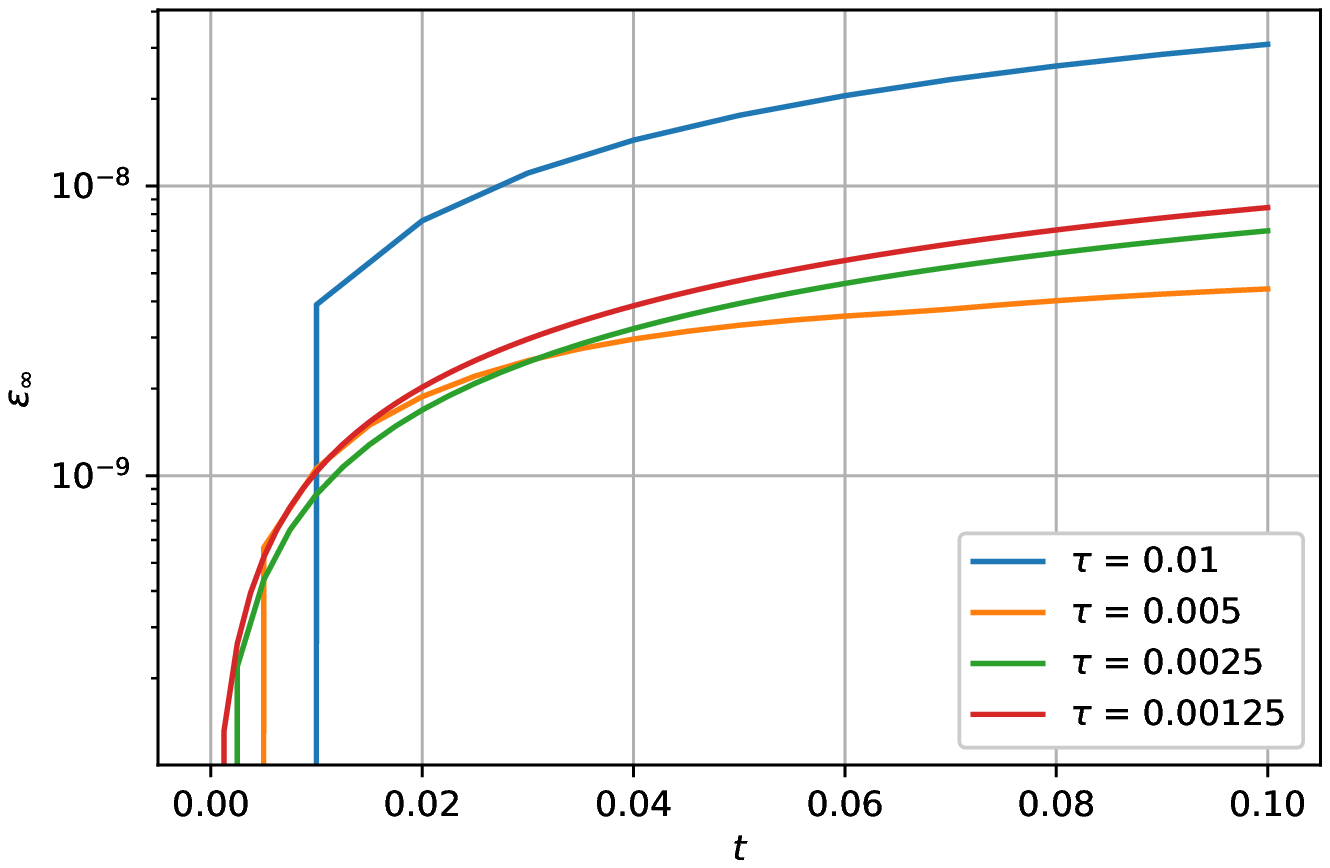}
\end{minipage}
\caption{Accuracy for the non-stationary problem: $\alpha = 0.25$, $m = 50$, $\sigma = 1$ (top) and $\sigma = 0.5$ (bottom).}
\label{f-4}
\end{figure}

\begin{figure}
\centering
\begin{minipage}{0.45\linewidth}
\centering
\includegraphics[width=\linewidth]{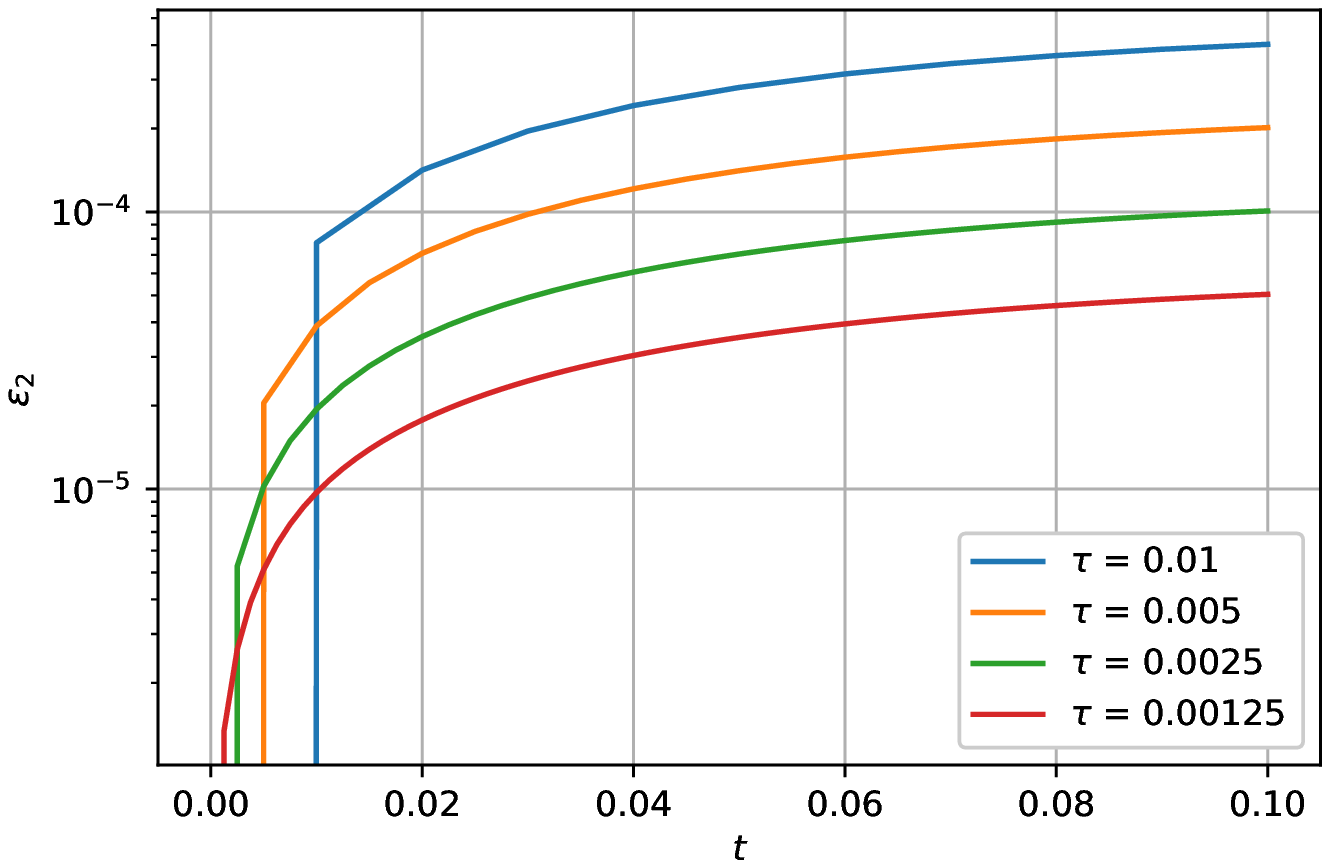}
\end{minipage}
\begin{minipage}{0.45\linewidth}
\centering
\includegraphics[width=\linewidth]{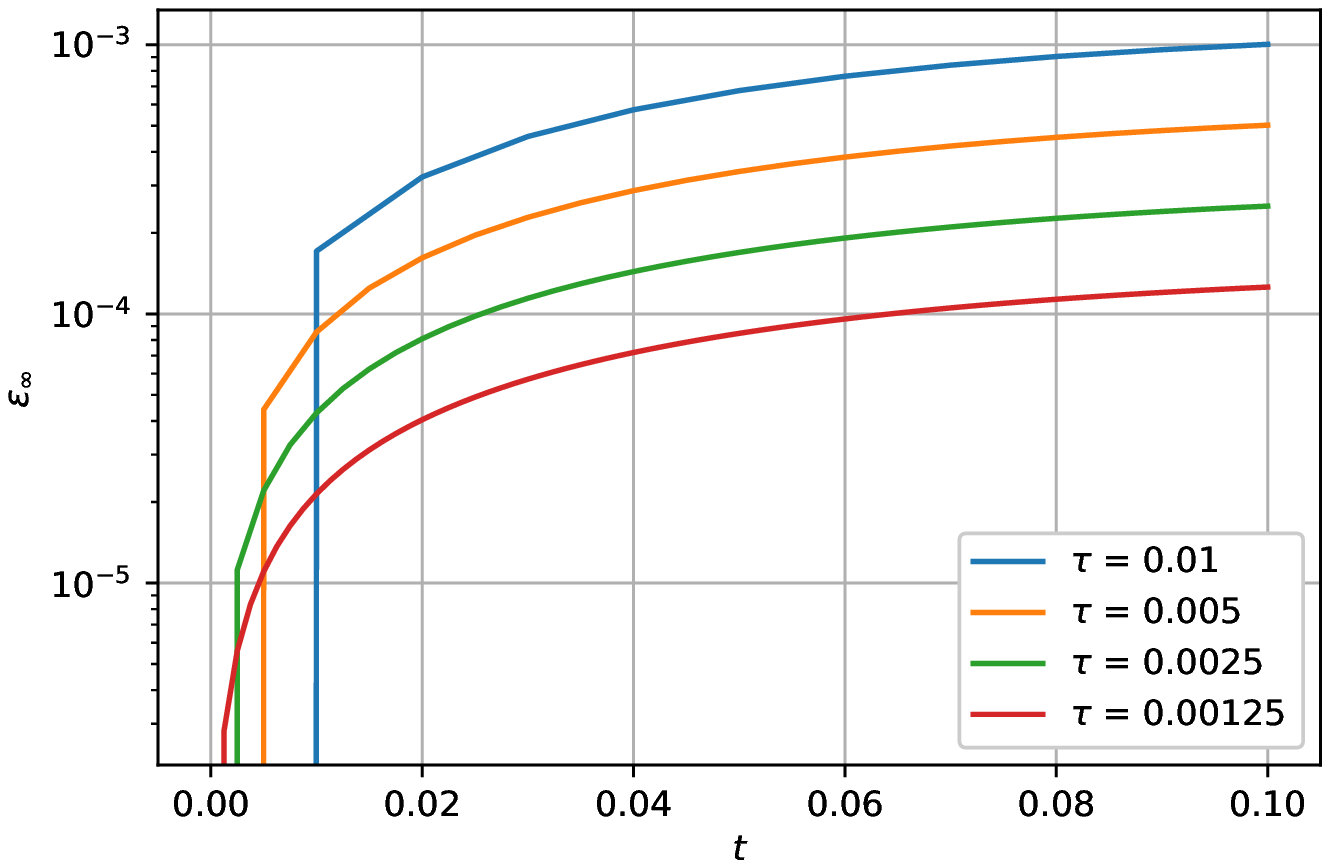}
\end{minipage}
\begin{minipage}{0.45\linewidth}
\centering
\includegraphics[width=\linewidth]{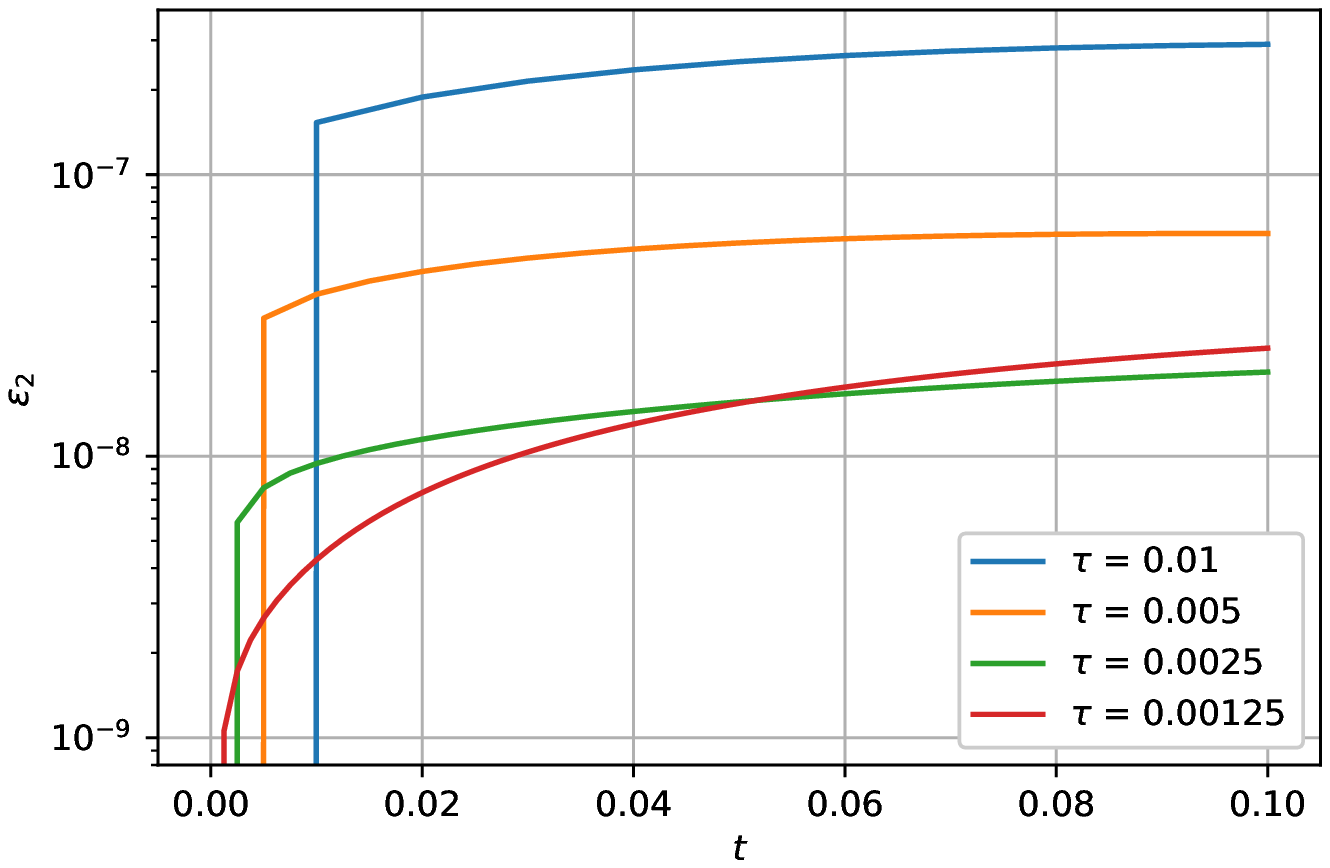}
\end{minipage}
\begin{minipage}{0.45\linewidth}
\centering
\includegraphics[width=\linewidth]{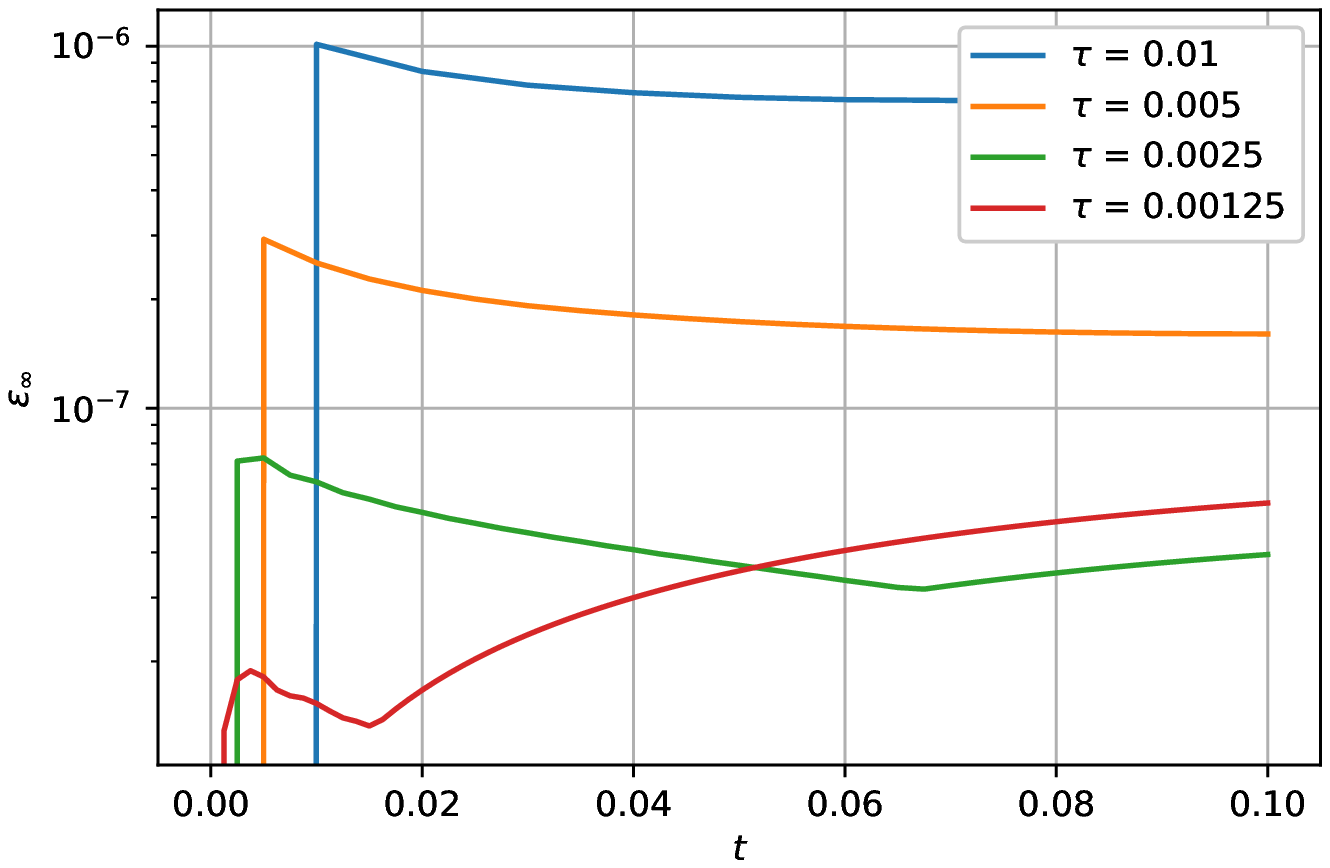}
\end{minipage}
\caption{Accuracy for the non-stationary problem: $\alpha = 0.5$, $m = 50$, $\sigma = 1$ (top) and $\sigma = 0.5$ (bottom).}
\label{f-5}
\end{figure}

\begin{figure}
\centering
\begin{minipage}{0.45\linewidth}
\centering
\includegraphics[width=\linewidth]{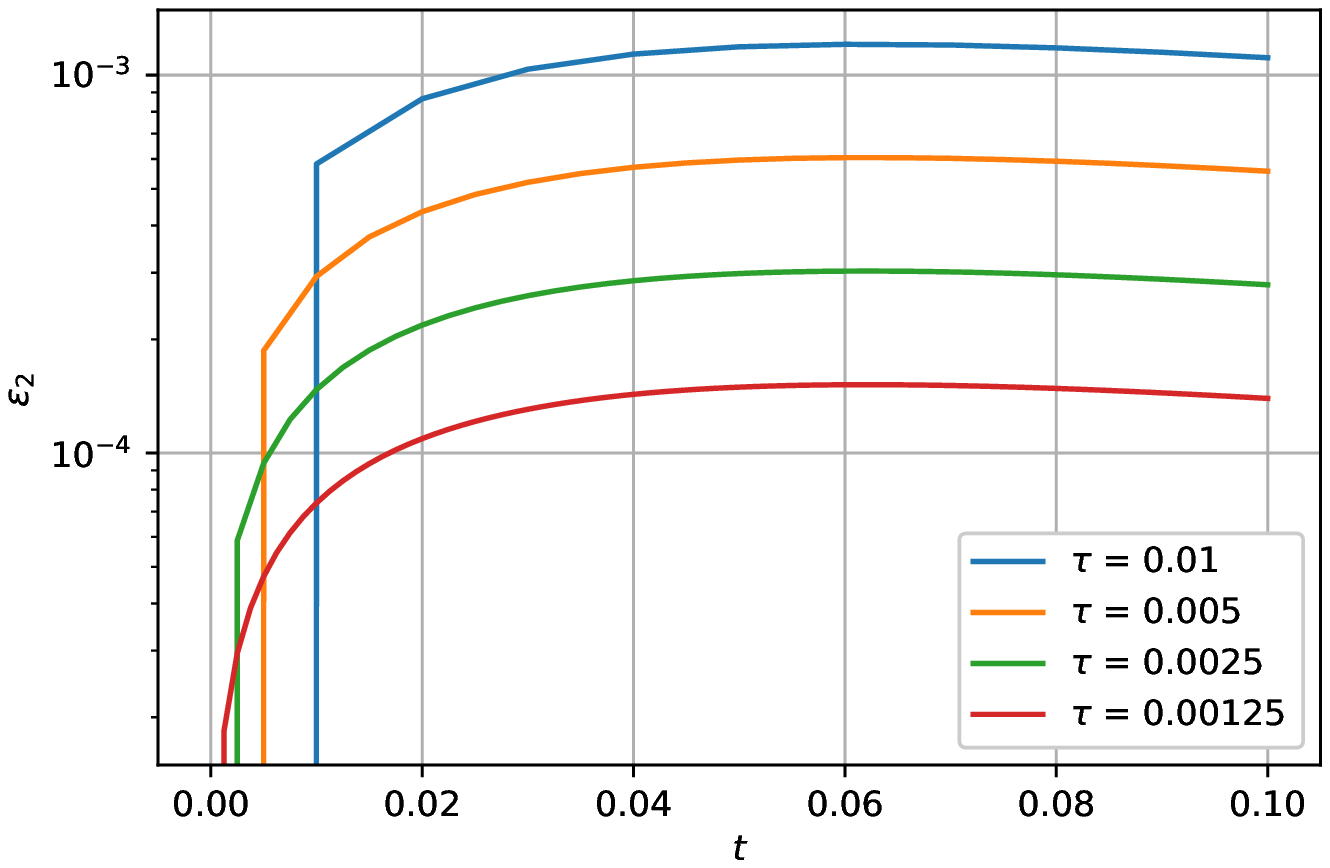}
\end{minipage}
\begin{minipage}{0.45\linewidth}
\centering
\includegraphics[width=\linewidth]{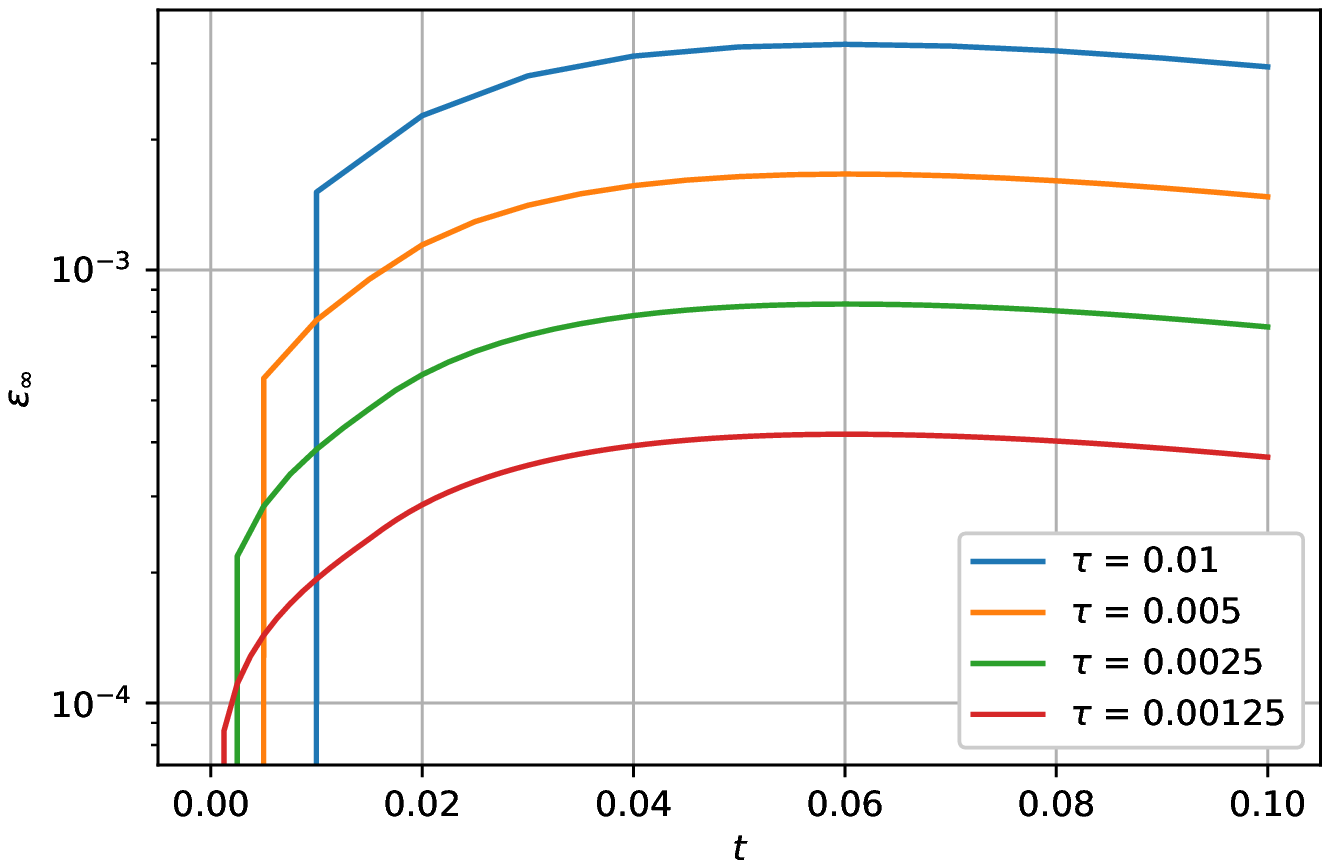}
\end{minipage}
\begin{minipage}{0.45\linewidth}
\centering
\includegraphics[width=\linewidth]{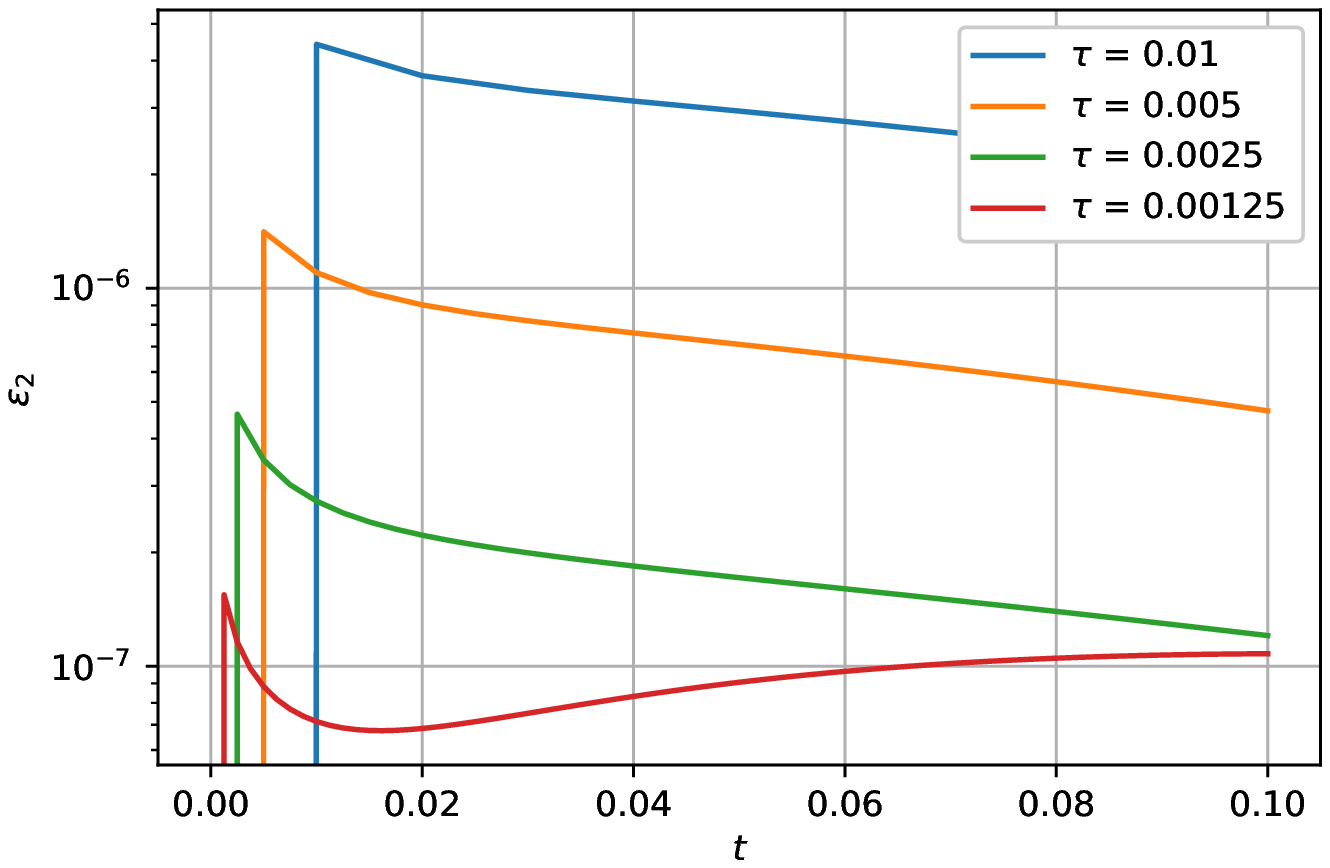}
\end{minipage}
\begin{minipage}{0.45\linewidth}
\centering
\includegraphics[width=\linewidth]{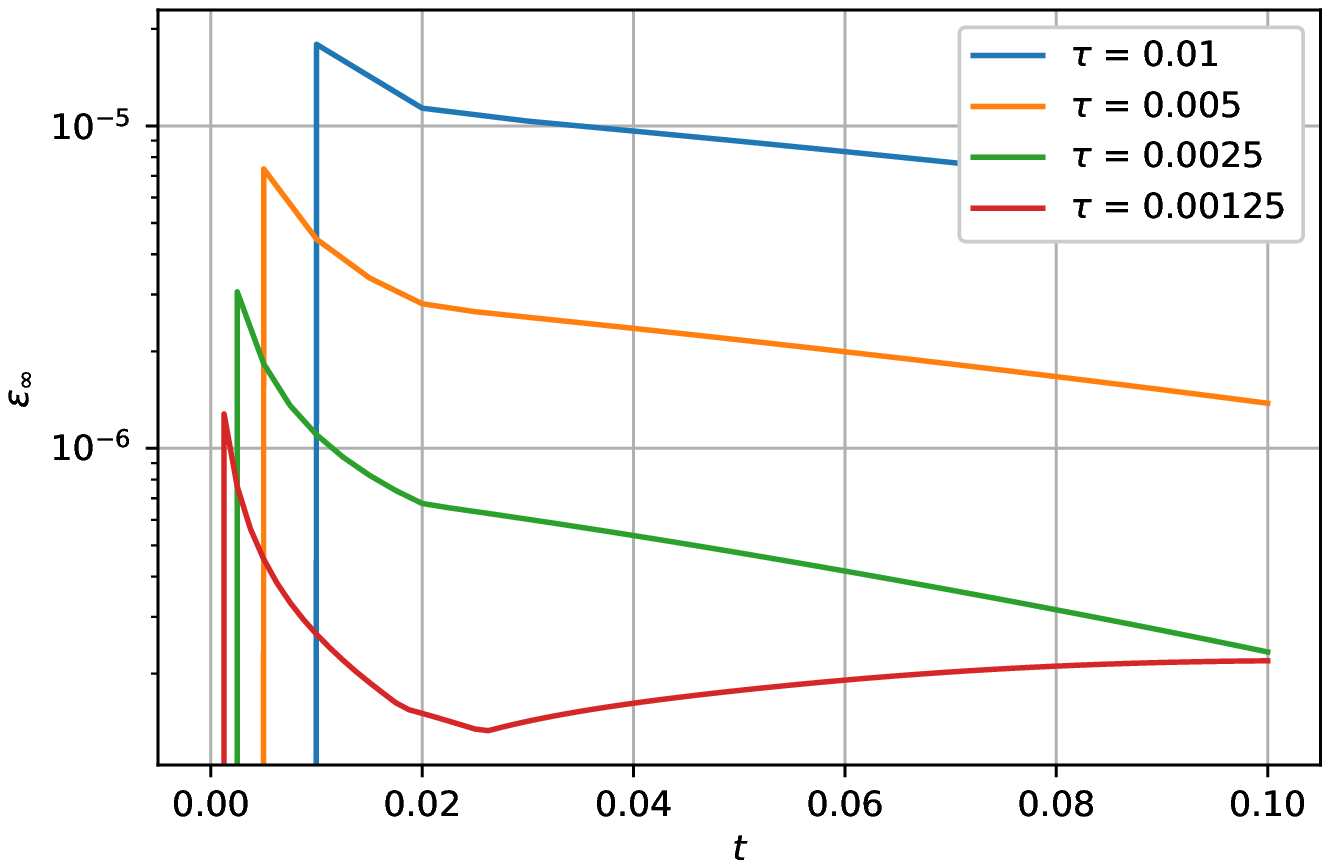}
\end{minipage}
\caption{Accuracy for the non-stationary problem: $\alpha = 0.75$, $m = 50$, $\sigma = 1$ (top) and $\sigma = 0.5$ (bottom).}
\label{f-6}
\end{figure}

\begin{figure}
\centering
\begin{minipage}{0.45\linewidth}
\centering
\includegraphics[width=\linewidth]{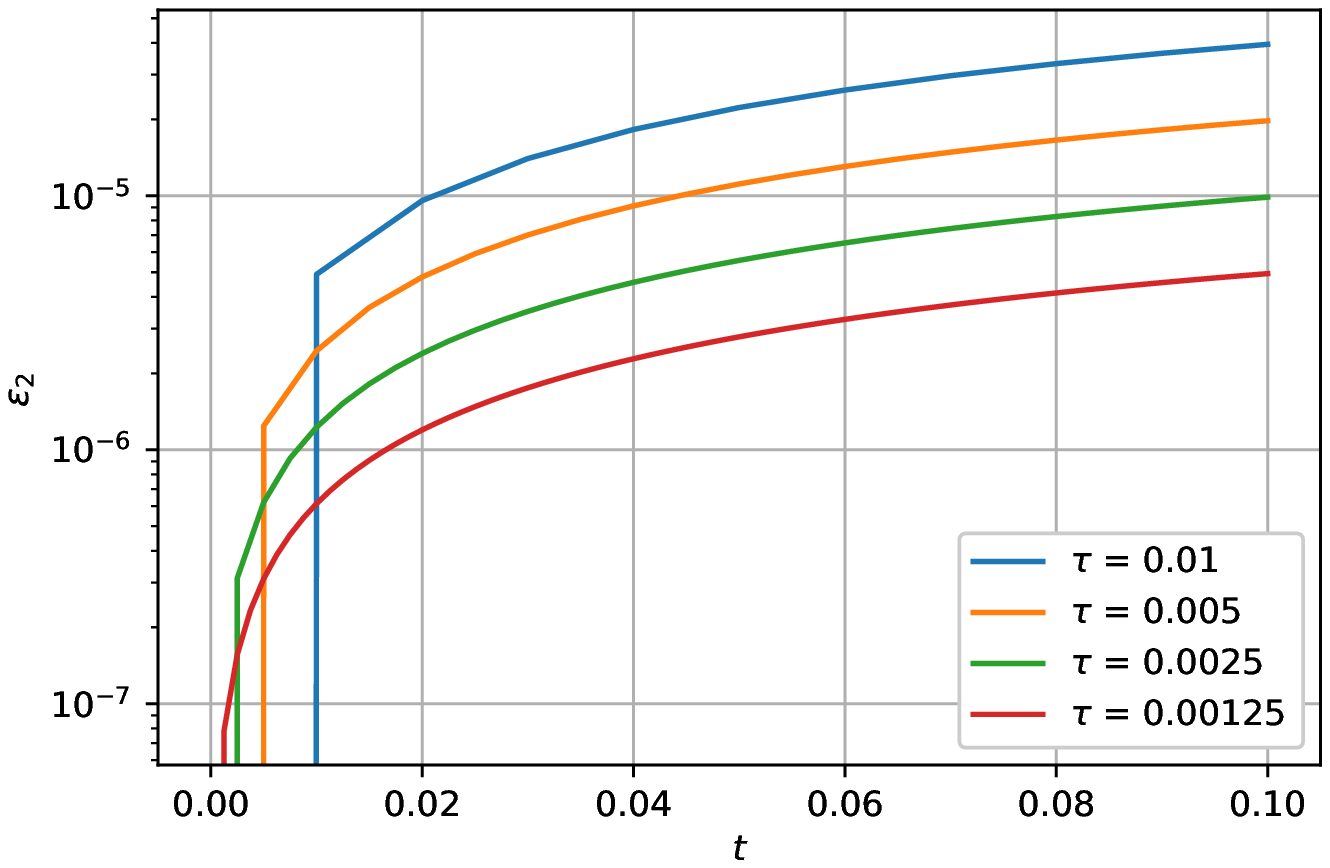}
\end{minipage}
\begin{minipage}{0.45\linewidth}
\centering
\includegraphics[width=\linewidth]{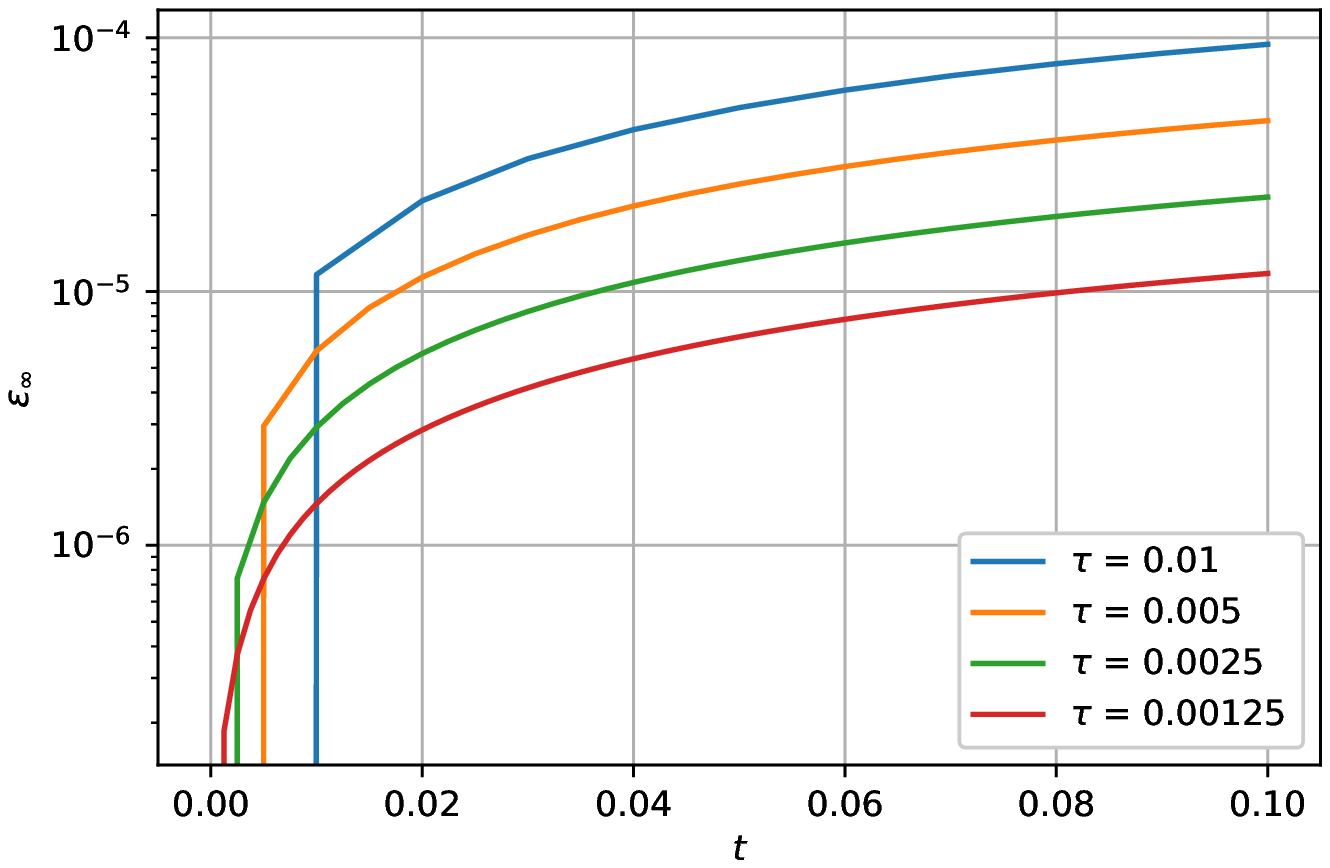}
\end{minipage}
\begin{minipage}{0.45\linewidth}
\centering
\includegraphics[width=\linewidth]{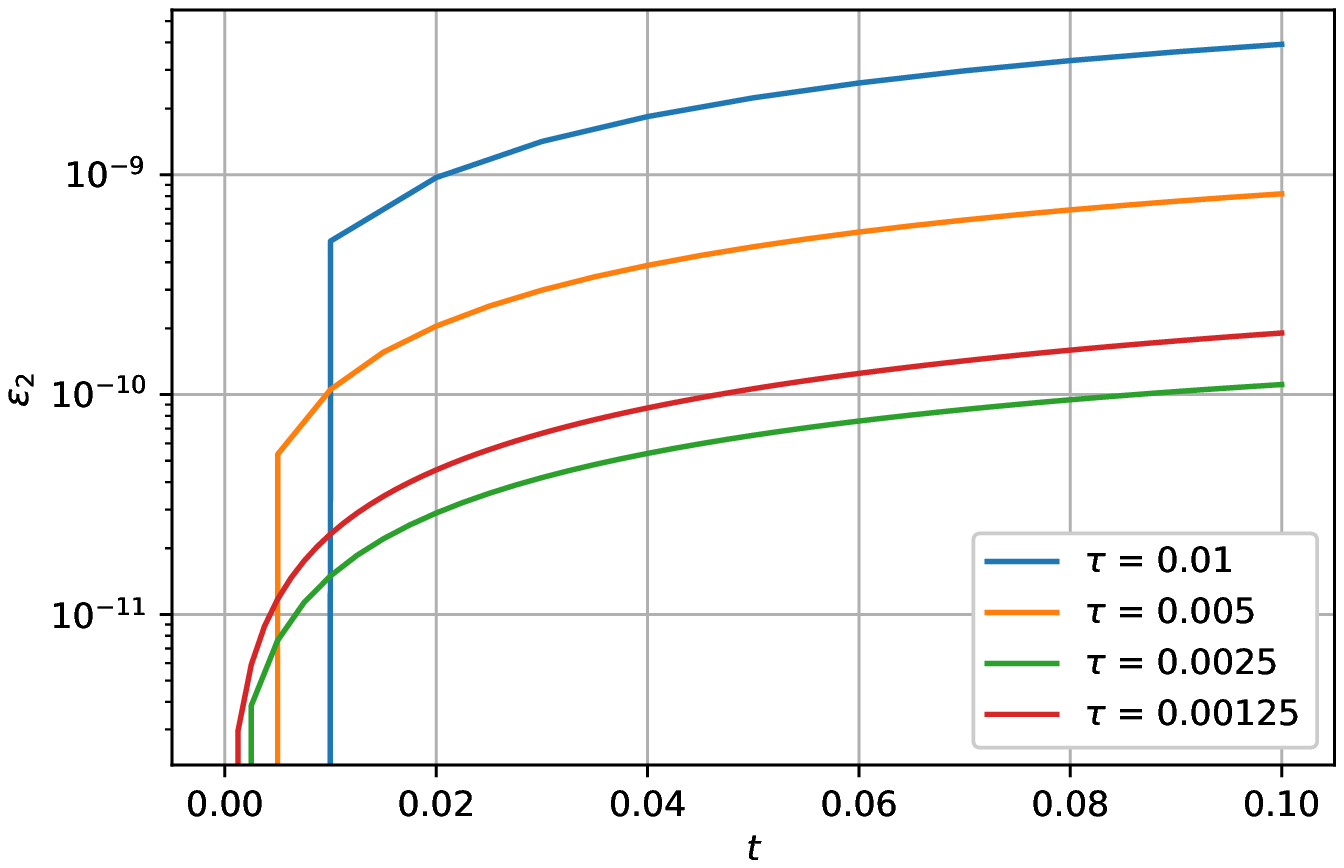}
\end{minipage}
\begin{minipage}{0.45\linewidth}
\centering
\includegraphics[width=\linewidth]{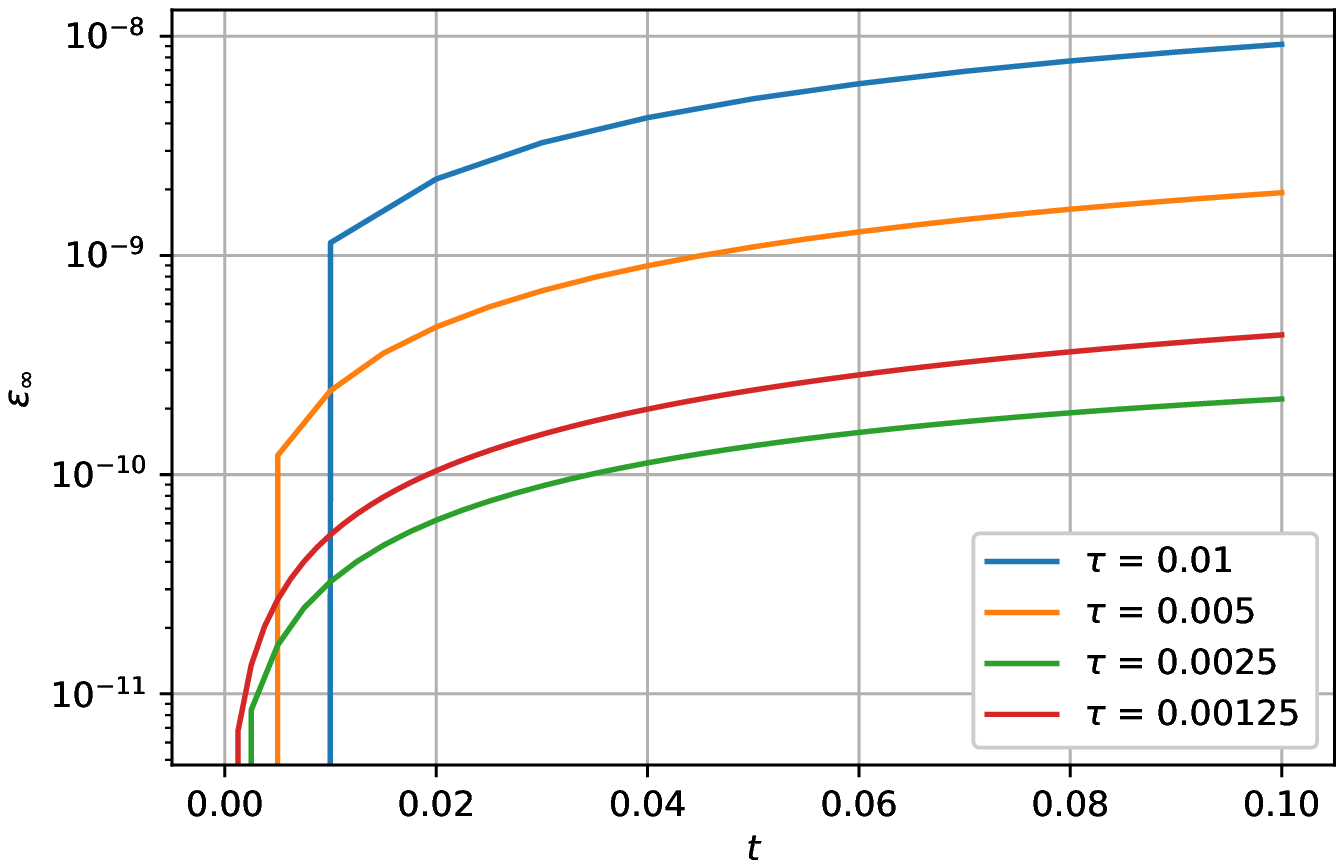}
\end{minipage}
\caption{Accuracy for the non-stationary problem: $\alpha = 0.25$, $m = 100$, $\sigma = 1$ (top) and $\sigma = 0.5$ (bottom).}
\label{f-7}
\end{figure}

The accuracy of the solution to the non-stationary problem is estimated by the absolute discrepancy at individual points in time:
\[
 \varepsilon_2(t^n) = \|u (\bm x,t^n) - w^n(\bm x)\|,
 \quad  \varepsilon_\infty (t^n)  = \|u (\bm x,t^n) - w^n(\bm x)\|_\infty ,
 \quad n = 0, \ldots, N . 
\] 
For the initial condition, we have
\[
 \|u^0 (\bm x)\| \approx 0.95238095 ,
 \quad \|u^0 (\bm x)\|_\infty \approx  2.19473708 . 
\] 
The dependence of the accuracy on the $\alpha$ parameter when using different time grids is shown in Figs.\ref{f-4}--\ref{f-6}. 
The calculations were performed using the Simpson quadrature formula when dividing into $m = 50$ intervals.
When using the scheme (\ref{20}), (\ref{21}) with $\sigma = 1$ convergence with the first order in $\tau$ is observed.
The splitting scheme with $\sigma = 0.5$ has a much higher accuracy.
In this case, the error in the rational approximation of the operator $A^{-\beta}$ is clearly manifested.
An illustration is Fig.\ref{f-7}, which shows the results of solving the problem with $\alpha = 0.25$ using a more accurate approximation with $m = 100$.

When solving multidimensional non-stationary problems for partial differential equations, we match the required accuracy of the approximate solution with discretization in space and time.
When considering problems with fractional powers of operators, it is necessary to take into account the error in approximating fractional powers of an operator.

\section{Conclusions} 

\begin{enumerate}
 \item We consider the Cauchy problem for a first-order evolution equation with a fractional power of the operator.
  Such nonlocal mathematical models are used to simulate phenomena and processes of various natures, and  they are actively discussed in the literature.
 \item The approach with rational approximation of the fractional power of an operator in various versions is widely used.
  In particular, in many works, the rational approximation is based on one or another integral representation of the fractional power operator.
  We construct an additive representation of the problem operator in the study of non-stationary problems.
 \item Two-level additive operator-difference schemes are proposed and investigated for a first-order evolution equation with a fractional power operator. When using rational approximation, the transition to a new level in time is realized by the sequential solution of standard evolutionary problems.
 \item Problems with an operator at the derivative in evolutionary equations of the first order are distinguished among more complex problems with a fractional power of the operator.
We construct unconditionally stable splitting schemes based on the perturbation of individual operator terms in rational approximation by a fractional power of the operator.
 \item We demonstrate the capabilities of the proposed splitting schemes for the numerical solution of a model two-dimensional problem in a rectangle with a fractional power of the Laplace operator.
We use conventional finite-difference approximations in space and two-level weighted schemes with time approximations.
We construct rational approximations based on the integral representation of the fractional power of the operator, which we proposed earlier.
\end{enumerate} 

\section*{Acknowledgements}

The publication has been prepared with support of the mega-grant of
the Russian Federation Government 14.Y26.31.0013 and the research grant 20-01-00207 of Russian Foundation of Basic Research.


\begin{thebibliography}{37}
\expandafter\ifx\csname natexlab\endcsname\relax\def\natexlab#1{#1}\fi
\providecommand{\url}[1]{\texttt{#1}}
\providecommand{\href}[2]{#2}
\providecommand{\path}[1]{#1}
\providecommand{\DOIprefix}{doi:}
\providecommand{\ArXivprefix}{arXiv:}
\providecommand{\URLprefix}{URL: }
\providecommand{\Pubmedprefix}{pmid:}
\providecommand{\doi}[1]{\href{http://dx.doi.org/#1}{\path{#1}}}
\providecommand{\Pubmed}[1]{\href{pmid:#1}{\path{#1}}}
\providecommand{\bibinfo}[2]{#2}
\ifx\xfnm\relax \def\xfnm[#1]{\unskip,\space#1}\fi
\bibitem[{Baleanu(2012)}]{baleanu2012fractional}
\bibinfo{author}{D.~Baleanu}, \bibinfo{title}{Fractional Calculus: Models and
  Numerical Methods}, \bibinfo{publisher}{World Scientific},
  \bibinfo{address}{New York}, \bibinfo{year}{2012}.
\bibitem[{Uchaikin(2013)}]{uchaikin}
\bibinfo{author}{V.~V. Uchaikin}, \bibinfo{title}{Fractional Derivatives for
  Physicists and Engineers}, \bibinfo{publisher}{Springer},
  \bibinfo{address}{Heidelberg}, \bibinfo{year}{2013}.
\bibitem[{Pozrikidis(2018)}]{Pozrikidis16}
\bibinfo{author}{C.~Pozrikidis}, \bibinfo{title}{The Fractional Laplacian},
  \bibinfo{publisher}{CRC Press}, \bibinfo{address}{Boca Raton},
  \bibinfo{year}{2018}.
\bibitem[{Knabner and Angermann(2003)}]{KnabnerAngermann2003}
\bibinfo{author}{P.~Knabner}, \bibinfo{author}{L.~Angermann},
  \bibinfo{title}{Numerical Methods for Elliptic and Parabolic Partial
  Differential Equations}, \bibinfo{publisher}{Springer}, \bibinfo{address}{New
  York}, \bibinfo{year}{2003}.
\bibitem[{Quarteroni and Valli(1994)}]{QuarteroniValli1994}
\bibinfo{author}{A.~Quarteroni}, \bibinfo{author}{A.~Valli},
  \bibinfo{title}{Numerical Approximation of Partial Differential Equations},
  \bibinfo{publisher}{Springer-Verlag}, \bibinfo{address}{Berlin},
  \bibinfo{year}{1994}.
\bibitem[{Birman and Solomjak(1987)}]{birman1987spectral}
\bibinfo{author}{M.~S. Birman}, \bibinfo{author}{M.~Z. Solomjak},
  \bibinfo{title}{Spectral theory of self-adjoint operators in Hilbert space},
  \bibinfo{publisher}{Kluwer academic publishers},
  \bibinfo{address}{Dordrecht}, \bibinfo{year}{1987}.
\bibitem[{Carracedo et~al.(2001)Carracedo, Alix, and
  Sanz}]{carracedo2001theory}
\bibinfo{author}{C.~M. Carracedo}, \bibinfo{author}{M.~S. Alix},
  \bibinfo{author}{M.~Sanz}, \bibinfo{title}{The Theory of Fractional Powers of
  Operators}, \bibinfo{publisher}{Elsevier}, \bibinfo{address}{Amsterdam},
  \bibinfo{year}{2001}.
\bibitem[{Higham(2008)}]{higham2008functions}
\bibinfo{author}{N.~J. Higham}, \bibinfo{title}{Functions of Matrices: Theory
  and Computation}, \bibinfo{publisher}{SIAM}, \bibinfo{address}{Philadelphia},
  \bibinfo{year}{2008}.
\bibitem[{Bonito et~al.(2018)Bonito, Borthagaray, Nochetto, Ot{\'a}rola, and
  Salgado}]{bonito2018numerical}
\bibinfo{author}{A.~Bonito}, \bibinfo{author}{J.~P. Borthagaray},
  \bibinfo{author}{R.~H. Nochetto}, \bibinfo{author}{E.~Ot{\'a}rola},
  \bibinfo{author}{A.~J. Salgado},
\newblock \bibinfo{title}{Numerical methods for fractional diffusion},
\newblock \bibinfo{journal}{Computing and Visualization in Science}
  \bibinfo{volume}{19} (\bibinfo{year}{2018}) \bibinfo{pages}{19--46}.
\bibitem[{Harizanov et~al.(2020)Harizanov, Lazarov, and
  Margenov}]{harizanov2020rev}
\bibinfo{author}{S.~Harizanov}, \bibinfo{author}{R.~Lazarov},
  \bibinfo{author}{S.~Margenov}, \bibinfo{title}{A Survey on Numerical Methods
  for Spectral Space-Fractional Diffusion Problems}, \bibinfo{type}{Technical
  Report} \bibinfo{number}{2010.02717}, arXiv, \bibinfo{year}{2020}.
\bibitem[{Stahl(2003)}]{stahl2003best}
\bibinfo{author}{H.~R. Stahl},
\newblock \bibinfo{title}{Best uniform rational approximation of $x^\alpha$ on
  [0, 1]},
\newblock \bibinfo{journal}{Acta Mathematica} \bibinfo{volume}{190}
  (\bibinfo{year}{2003}) \bibinfo{pages}{241--306}.
\bibitem[{Harizanov et~al.(2020{\natexlab{a}})Harizanov, Lazarov, Margenov, and
  Marinov}]{harizanov2020num}
\bibinfo{author}{S.~Harizanov}, \bibinfo{author}{R.~Lazarov},
  \bibinfo{author}{S.~Margenov}, \bibinfo{author}{P.~Marinov},
\newblock \bibinfo{title}{Numerical solution of fractional diffusion–reaction
  problems based on {BURA}},
\newblock \bibinfo{journal}{Computers \& Mathematics with Applications}
  \bibinfo{volume}{80} (\bibinfo{year}{2020}{\natexlab{a}})
  \bibinfo{pages}{316--331}.
\bibitem[{Harizanov et~al.(2020{\natexlab{b}})Harizanov, Lazarov, Marinov,
  Margenov, and Pasciak}]{harizanov2020anal}
\bibinfo{author}{S.~Harizanov}, \bibinfo{author}{R.~Lazarov},
  \bibinfo{author}{P.~Marinov}, \bibinfo{author}{S.~Margenov},
  \bibinfo{author}{J.~Pasciak},
\newblock \bibinfo{title}{Analysis of numerical methods for spectral fractional
  elliptic equations based on the best uniform rational approximation},
\newblock \bibinfo{journal}{Journal of Computational Physics}
  \bibinfo{volume}{408} (\bibinfo{year}{2020}{\natexlab{b}})
  \bibinfo{pages}{109285}.
\bibitem[{Frommer et~al.(2014)Frommer, G\"{u}ttel, and
  Schweitzer}]{frommer2014efficient}
\bibinfo{author}{A.~Frommer}, \bibinfo{author}{S.~G\"{u}ttel},
  \bibinfo{author}{M.~Schweitzer},
\newblock \bibinfo{title}{Efficient and stable {A}rnoldi restarts for matrix
  functions based on quadrature},
\newblock \bibinfo{journal}{SIAM Journal on Matrix Analysis and Applications}
  \bibinfo{volume}{35} (\bibinfo{year}{2014}) \bibinfo{pages}{661--683}.
\bibitem[{Bonito and Pasciak(2015)}]{bonito2015numerical}
\bibinfo{author}{A.~Bonito}, \bibinfo{author}{J.~Pasciak},
\newblock \bibinfo{title}{Numerical approximation of fractional powers of
  elliptic operators},
\newblock \bibinfo{journal}{Mathematics of Computation} \bibinfo{volume}{84}
  (\bibinfo{year}{2015}) \bibinfo{pages}{2083--2110}.
\bibitem[{Aceto and Novati(2019)}]{Aceto2019}
\bibinfo{author}{L.~Aceto}, \bibinfo{author}{P.~Novati},
\newblock \bibinfo{title}{Rational approximations to fractional powers of
  self-adjoint positive operators},
\newblock \bibinfo{journal}{Numerische Mathematik} \bibinfo{volume}{143}
  (\bibinfo{year}{2019}) \bibinfo{pages}{1--16}.
\bibitem[{Balakrishnan(1960)}]{balakrishnan1960fractional}
\bibinfo{author}{A.~V. Balakrishnan},
\newblock \bibinfo{title}{Fractional powers of closed operators and the
  semigroups generated by them},
\newblock \bibinfo{journal}{Pacific Journal of Mathematics}
  \bibinfo{volume}{10} (\bibinfo{year}{1960}) \bibinfo{pages}{419--437}.
\bibitem[{Vabishchevich(2020)}]{vabishchevich2020}
\bibinfo{author}{P.~N. Vabishchevich},
\newblock \bibinfo{title}{Approximation of a fractional power of an elliptic
  operator},
\newblock \bibinfo{journal}{Linear Algebra and its Applications}
  \bibinfo{volume}{27} (\bibinfo{year}{2020}) \bibinfo{pages}{e2287}.
\bibitem[{Nochetto et~al.(2015)Nochetto, Ot{\'a}rola, and
  Salgado}]{nochetto2015pde}
\bibinfo{author}{R.~H. Nochetto}, \bibinfo{author}{E.~Ot{\'a}rola},
  \bibinfo{author}{A.~J. Salgado},
\newblock \bibinfo{title}{A {PDE} approach to fractional diffusion in general
  domains: a priori error analysis},
\newblock \bibinfo{journal}{Foundations of Computational Mathematics}
  \bibinfo{volume}{15} (\bibinfo{year}{2015}) \bibinfo{pages}{733--791}.
\bibitem[{Nochetto et~al.(2016)Nochetto, Otarola, and
  Salgado}]{nochetto2016pde}
\bibinfo{author}{R.~H. Nochetto}, \bibinfo{author}{E.~Otarola},
  \bibinfo{author}{A.~J. Salgado},
\newblock \bibinfo{title}{A {PDE} approach to space-time fractional parabolic
  problems},
\newblock \bibinfo{journal}{SIAM Journal on Numerical Analysis}
  \bibinfo{volume}{54} (\bibinfo{year}{2016}) \bibinfo{pages}{848--873}.
\bibitem[{Caffarelli and Silvestre(2007)}]{Caffarelli}
\bibinfo{author}{L.~Caffarelli}, \bibinfo{author}{L.~Silvestre},
\newblock \bibinfo{title}{An extension problem related to the fractional
  {L}aplacian},
\newblock \bibinfo{journal}{Communications in Partial Differential Equations}
  \bibinfo{volume}{32} (\bibinfo{year}{2007}) \bibinfo{pages}{1245--1260}.
\bibitem[{Vabishchevich(2015)}]{vabishchevich2014numerical}
\bibinfo{author}{P.~N. Vabishchevich},
\newblock \bibinfo{title}{Numerically solving an equation for fractional powers
  of elliptic operators},
\newblock \bibinfo{journal}{Journal of Computational Physics}
  \bibinfo{volume}{282} (\bibinfo{year}{2015}) \bibinfo{pages}{289--302}.
\bibitem[{Duan et~al.(2019)Duan, Lazarov, and Pasciak}]{duan2018numerical}
\bibinfo{author}{B.~Duan}, \bibinfo{author}{R.~Lazarov},
  \bibinfo{author}{J.~Pasciak},
\newblock \bibinfo{title}{Numerical approximation of fractional powers of
  elliptic operators},
\newblock \bibinfo{journal}{IMA J. Numerical Analysis} \bibinfo{volume}{40}
  (\bibinfo{year}{2019}) \bibinfo{pages}{1746--1771}.
\bibitem[{\v{C}iegis and Vabishchevich(2020)}]{CegVab2019}
\bibinfo{author}{R.~\v{C}iegis}, \bibinfo{author}{P.~N. Vabishchevich},
\newblock \bibinfo{title}{Two-level schemes of cauchy problem method for
  solving fractional powers of elliptic operators},
\newblock \bibinfo{journal}{Computers \& Mathematics with Applications}
  \bibinfo{volume}{80} (\bibinfo{year}{2020}) \bibinfo{pages}{305--315}.
\bibitem[{Hofreither(2020)}]{Hofreither2020}
\bibinfo{author}{C.~Hofreither},
\newblock \bibinfo{title}{A unified view of some numerical methods for
  fractional diffusion},
\newblock \bibinfo{journal}{Computers \& Mathematics with Applications}
  \bibinfo{volume}{80} (\bibinfo{year}{2020}) \bibinfo{pages}{332--350}.
\bibitem[{Yagi(2009)}]{yagi2009abstract}
\bibinfo{author}{A.~Yagi}, \bibinfo{title}{Abstract Parabolic Evolution
  Equations and Their Applications}, \bibinfo{publisher}{Springer},
  \bibinfo{address}{Berlin}, \bibinfo{year}{2009}.
\bibitem[{Samarskii(2001)}]{SamarskiiTheory}
\bibinfo{author}{A.~A. Samarskii}, \bibinfo{title}{The Theory of Difference
  Schemes}, \bibinfo{publisher}{Marcel Dekker}, \bibinfo{address}{New York},
  \bibinfo{year}{2001}.
\bibitem[{Harizanov et~al.(2020)Harizanov, Lazarov, Margenov, and
  Marinov}]{harizanov2020reac}
\bibinfo{author}{S.~Harizanov}, \bibinfo{author}{R.~Lazarov},
  \bibinfo{author}{S.~Margenov}, \bibinfo{author}{P.~Marinov},
\newblock \bibinfo{title}{Numerical solution of fractional diffusion–reaction
  problems based on bura},
\newblock \bibinfo{journal}{Computers \& Mathematics with Applications}
  \bibinfo{volume}{80} (\bibinfo{year}{2020}) \bibinfo{pages}{316–331}.
\bibitem[{Aceto and Novati(2017)}]{Aceto2017}
\bibinfo{author}{L.~Aceto}, \bibinfo{author}{P.~Novati},
\newblock \bibinfo{title}{Rational approximation to the fractional {L}aplacian
  operator in reaction-diffusion problems},
\newblock \bibinfo{journal}{SIAM Journal on Scientific Computing}
  \bibinfo{volume}{39} (\bibinfo{year}{2017}) \bibinfo{pages}{A214--A228}.
\bibitem[{Vabishchevich(2018)}]{Vabishchevich2018}
\bibinfo{author}{P.~N. Vabishchevich},
\newblock \bibinfo{title}{Numerical solution of time-dependent problems with
  fractional power elliptic operator},
\newblock \bibinfo{journal}{Computational Methods in Applied Mathematics}
  \bibinfo{volume}{18} (\bibinfo{year}{2018}) \bibinfo{pages}{111--128}.
\bibitem[{Vabishchevich(2016)}]{Vabishchevich2016}
\bibinfo{author}{P.~N. Vabishchevich},
\newblock \bibinfo{title}{Numerical solution of non-stationary problems for a
  space-fractional diffusion equation},
\newblock \bibinfo{journal}{Fractional Calculus and Applied Analysis}
  \bibinfo{volume}{19} (\bibinfo{year}{2016}) \bibinfo{pages}{116--139}.
\bibitem[{Samarskii et~al.(2002)Samarskii, Matus, and
  Vabishchevich}]{SamarskiiMatusVabischevich2002}
\bibinfo{author}{A.~A. Samarskii}, \bibinfo{author}{P.~P. Matus},
  \bibinfo{author}{P.~N. Vabishchevich}, \bibinfo{title}{Difference Schemes
  with Operator Factors}, \bibinfo{publisher}{Kluwer Academic},
  \bibinfo{address}{Dordrecht}, \bibinfo{year}{2002}.
\bibitem[{Vabishchevich(2018)}]{Vabishchevich2018a}
\bibinfo{author}{P.~N. Vabishchevich},
\newblock \bibinfo{title}{Numerical solution of time-dependent problems with a
  fractional-power elliptic operator},
\newblock \bibinfo{journal}{Computational Mathematics and Mathematical Physics}
  \bibinfo{volume}{58} (\bibinfo{year}{2018}) \bibinfo{pages}{394--409}.
\bibitem[{Marchuk(1990)}]{Marchuk1990}
\bibinfo{author}{G.~I. Marchuk},
\newblock \bibinfo{title}{Splitting and alternating direction methods},
\newblock in: \bibinfo{editor}{P.~G. Ciarlet}, \bibinfo{editor}{J.-L. Lions}
  (Eds.), \bibinfo{booktitle}{Handbook of Numerical Analysis, Vol. I},
  \bibinfo{publisher}{North-Holland}, \bibinfo{year}{1990}, pp.
  \bibinfo{pages}{197--462}.
\bibitem[{Vabishchevich(2013)}]{VabishchevichAdditive}
\bibinfo{author}{P.~N. Vabishchevich}, \bibinfo{title}{Additive
  Operator-Difference Schemes: Splitting Schemes}, \bibinfo{publisher}{Walter
  de Gruyter GmbH}, \bibinfo{address}{Berlin, Boston}, \bibinfo{year}{2013}.
\bibitem[{Samarskii and Nikolaev(1989)}]{SamarskiiNikolaev1978}
\bibinfo{author}{A.~A. Samarskii}, \bibinfo{author}{E.~S. Nikolaev},
  \bibinfo{title}{Numerical methods for grid equations. Vol. I, II},
  \bibinfo{publisher}{Birkhauser Verlag}, \bibinfo{address}{Basel},
  \bibinfo{year}{1989}.
\bibitem[{Ralston and Rabinowitz(2001)}]{Rabinowitz}
\bibinfo{author}{A.~Ralston}, \bibinfo{author}{P.~Rabinowitz},
  \bibinfo{title}{A First Course in Numerical Analysis},
  \bibinfo{publisher}{Dover Publications}, \bibinfo{address}{Mineola, NY},
  \bibinfo{year}{2001}.

\end{thebibliography}
\end{document}